\documentclass[review]{elsarticle}
\usepackage{lineno,hyperref}
\modulolinenumbers[5]
\usepackage{amsmath, amsthm}
\usepackage{amssymb}
\usepackage{graphics}
\usepackage{color}
\input{epsf.sty}
\usepackage{epsfig}
\newtheorem{Definition_}{Definition}

\newtheorem{Theorem_}{Theorem}

\newtheorem{Property_}{Property}
\newtheorem{Remark_}{Remark}
\newtheorem{Lemma_}{Lemma}

\newtheorem{Corollary_}{Corollary}
\newtheorem{Assumption_}{Assumption}

\newcommand{\commentout}[1]{}
\newcommand{\be}{\begin{equation}}
\newcommand{\ee}{\end{equation}}
\newcommand{\bea}{\begin{eqnarray}}
\newcommand{\eea}{\end{eqnarray}}
\newcommand{\bean}{\begin{eqnarray*}}
\newcommand{\eean}{\end{eqnarray*}}

\newenvironment{pf}{\vspace{0.1in}\noindent \textbf{Proof $\!\!\!$} ~}{~
 \hspace*{\fill} $\square$ }
\newcommand{\ok}[1]{{\color{black}#1}}
\let\textdoublequote=\"
\let\mathdoublequote=\ddot
\DeclareRobustCommand{\"}[1]{\ifmmode {\mathdoublequote{#1}}\else  \textdoublequote{#1}\fi}
\let\mathddot\ddot
\DeclareRobustCommand{\ddot}[1]{\ifmmode{\mathddot{#1}}\else
                                        \textdoublequote{#1}\fi}
\DeclareRobustCommand{\bold}[1]{\ifmmode \mathbf{#1}\else\textbf{#1}\fi}
\DeclareRobustCommand{\italic}[1]{\ifmmode\mathit{#1}\else\textit{#1}\fi}
\DeclareRobustCommand{\roman}[1]{\ifmmode\mathrm{#1}\else\textrm{#1}\fi}
\DeclareMathAlphabet{\mathbi}{\encodingdefault}{\rmdefault}{\bfdefault}{\itdefault}
\DeclareRobustCommand{\bit}[1]{\ifmmode\mathbi{#1}\else\textbf{\textit{#1}}\fi}

\long\def\@makecaption#1#2{%
  \vskip\abovecaptionskip
  \sbox\@tempboxa{#1: #2}%
  \ifdim \wd\@tempboxa >\hsize #1: #2\par \else \global
  \@minipagefalse
    \hb@xt@\hsize{\box\@tempboxa\hfil}%
  \fi \vskip\belowcaptionskip}

\newcommand{\lchi}[1]{\mbox{\large $\chi$}}

\begin{document}

\begin{frontmatter}
\title{\ok{Ameso Optimization: a Relaxation of Discrete Midpoint Convexity}}

\author{Wen Chen}
\address{Providence Business School, Providence College, Providence, RI 02908}
\ead{wchen@providence.edu}

\author{Odysseas Kanavetas}
\address{Department of Management Science and Information Systems,
Rutgers Business School,  Newark and New Brunswick, 180 University
Avenue, Newark, NJ  07102-1895} \ead{okanavetas@business.rutgers.edu}

\author{Michael N. Katehakis}
\address{Department of Management Science and Information Systems,
Rutgers Business School,  Newark and New Brunswick, 180 University
Avenue, Newark, NJ  07102-1895} \ead{mnk@andromeda.rutgers.edu}

\date{}
%\subjclass{Primary 60J10, 28A33; Secondary 28C15}

\begin{keyword}
Discrete Optimization, Integral Optimization, Recursive
Procedure
\end{keyword}

\begin{abstract}
In this paper we introduce the {\sl Ameso} optimization problem, a special class of discrete optimization
problems.  We establish its basic properties and   investigate the relation between {\sl Ameso} optimization and the convex optimization.  Further, we design an algorithm to solve  
a  multi-dimensional {\sl Ameso} problem by solving  a sequence of one-dimensional {\sl Ameso} problems. Finally, we  demonstrate  how the  knapsack problem can be solved using the 
{\sl Ameso} optimization framework. 
% ii) we defined the conditional pair of an  Ameso($C$) pair and we
%showed that it is also an Ameso($C$) pair,

\end{abstract}
\end{frontmatter}
%\linenumbers
\section{Introduction }
In this paper
we introduce a new class of discrete optimization problems the {\sl
Ameso} optimization problems. For the one dimensional  case we    
show that  any optimal point can be determined by simple to
verify, optimality conditions. Furthermore, we  construct the
Ameso Recursive Procedure ({\bf ARP}) that solves   Ameso
optimization problems without necessarily performing complete
enumeration. Parallel implementations of the {\bf ARP} can easily be
done, c.f. \cite{BDP}.  Since this is a new class of problems there
is no directly related literature. However, since an Ameso problem can be a generalization and relaxation of midpoint convexity there are algorithms proposed in other papers that employ the proximity framework while using descent algorithms for discrete midpoint convex functions, c.f. \cite{MMTT}. In this paper, they are trying to highlight discrete midpoint convexity as a
unifying framework for convexity concepts of functions on the integer lattice $\mathbb{Z}^n$, and to investigate structural and algorithmic properties of functions defined by versions of discrete midpoint convexity. Using these structural properties they develop a proximity-scaling based algorithm for the minimization of locally and globally discrete midpoint convex functions. Scaling and proximity algorithms are successful for discrete optimization problems such as resource allocation problems \cite{H}, \cite{HS}, \cite{IK}, \cite{KSI} and convex network flow problems \cite{AMO}, \cite{IMM}, \cite{IS}.

Also, based on the midpoint convexity there are many other approaches to nonlinear integer optimization as in \cite{LL},  and other more algebraic methods have been developed in the last two decades, as described in \cite{DEL}, \cite{HKW}, \cite{O}. Finally, one can find a lot of applications where proving that a model is Ameso we can have very simple algorithm to obtain the optimal solution, as described in \cite{Z}. This paper considers a one-period assemble-to-order system with stochastic demand and uses combinatorial optimization and discrete convexity, to decrease the computational complexity, for certain specially structured models.

%***********
We provide an illustration of an  application of the Ameso optimization   
using a version of the knapsack problem cf.   \cite{GilmoreGomony1961}, ~\cite{GilmoreGomony1963},~\cite{GilmoreGomony1966}. This is a 
classical discrete optimization problem, and solution techniques for it include 
the branch and bound algorithm ~\cite{Kolesar1967}, the renewal algorithm ~\cite{ShapiroWanger1967}, dynamic programming ~\cite{RossTsang1989}, ~\cite{Denardo1982} etc.
Since many discrete optimization problems can be reduced to knapsack problem, we thus demonstrate 
the  possibility  of using the Ameso based algorithms described herein, as an additional tool for 
such problems. 

% Hence, our Ameso function may provide another consideration of such a class of discrete optimization. The contribution of the paper connecting Ameso function with revised knapsack problem is to discuss the possibility application of ameso function in the discrete optimization area in the future study.
%********

The rest of the paper is organized as follows. In section 2, we
define the Ameso optimization problem and we discuss its properties.
In section 2.1 there is the relationship between convexity, midpoint convexity and Ameso optimization problem. Section 2.2 is devoted to the one dimensional case and in section
2.3 we discuss the main property of the high dimensional case and we
use this property to design a procedure to solve Ameso(C)
optimization problems.  Also, we present two examples that illustrate the
performance of the proposed procedure. In Section 3, we discuss the knapsack formulation.
Finally, in the conclusions we discuss the relaxation that  Ameso(1)  provides to the  midpoint convexity and potential benefits. 

\section{The Ameso($C$) Optimization Problem}

Given a subset $D^n $ of the $n$ - dimensional integers,
$D^n\subseteq\mathbb{Z}^n,$ and a real function $f$ defined on
$D^n,$ we define the following:

\begin{Definition_}\label{def_Amesoset}
%\begin{enumerate}
  %\item
%\textbf{Ameso set:}

\ok{$D^n$ is called {\bf Ameso set}, if it satisfies the following condition
\bea
\left\lceil\frac{\vec{x}+\vec{y}}{2}\right\rceil,\left\lfloor\frac{\vec{x}+\vec{y}}{2}\right\rfloor\in
D^n \mbox{,  for all $\vec{x},\vec{y}\in
D^n$.} \label{eqn_def_D}
\eea
}
\end{Definition_}

\begin{Definition_}\label{def_Amesopair}
\ok{$(D^n,f)$ is called an {\bf Ameso($C$) pair}, if and
only if it satisfies the following conditions}

\begin{itemize}
\item[-] the domain $D^n$ of the function $f(\cdot)$  is
an Ameso set,

\item[-]  $f(\cdot)$ has a lower  bound and

\item[-] there exists $C\geq 0$ such that the following holds for all $\vec{x},\vec{y}\in D^n$
\bea
f(\vec{x})+f(\vec{y})+C\geq
f(\lceil\frac{\vec{x}+\vec{y}}{2}\rceil)+f(\lfloor\frac{\vec{x}+\vec{y}}{2}\rfloor)
.\label{eqn_def_AmesoPair}
\eea

\end{itemize}
\end{Definition_}

 \begin{Definition_}\label{def_Amesooptimization}
\ok{Minimization of $f(\vec{x})$ subject to $\vec{x}\in D^n$
 is an {\bf Ameso($C$) optimization} problem, if  $(D^n,f) $ is an
Ameso($C$) pair.}
\end{Definition_}

{\bf Notation.} For notational simplicity in the sequel for any integers $a$ and $b$  we will use notation $[a,b]$, $[a,b)$, and $(a,b]$  to denote respectively the sets of integers:
$\{a,a+1,\ldots,b\}$, $\{a,a+1,\ldots,b-1\}$, and  $\{a+1,\ldots,b\} .$

For better understanding of the definition of an Ameso set we provide some examples below.
Also, note that for notational simplicity we use $D^n$  to denote sets
that are subsets of  $\mathbb{Z}^n$  that are not necessarily
products of identical subsets of $\mathbb{Z} $ c.f. 
 the Example 1 below.

 \textbf{Example 1:} Given any integers $a_i,b_i$ where $a_i<b_i$, if we let
$D_i=[a_i,b_i] ,$ $i=1,...,n$, then $D^n=D_1\times
\cdots \times D_n$ is an Ameso set.

 \textbf{Example 2:} The sets\\
  $A_1=\{(x_1,x_2)\, : \, x_2-x_1\geq 3, x_2\leq 10, \, x_1,x_2\in
 \mathbb{Z}^+\}$,\\
 $A_2=\{(x_1,x_2)\, : \, x_2-x_1\geq 0, x_2\leq 10, x_1\leq 5 \, x_1,x_2\in
 \mathbb{Z}^+\}$, \\
 $A_3=\{(x_1,x_2)\, : \, x_2-x_1\leq 10, \, x_1,x_2\in
 \mathbb{Z}^+\}$, \\
 are Ameso sets. However, the sets\\
 $A_4=\{(x_1,x_2)\, : \, 3x_2-x_1\geq 0,   x_1,x_2\in
 \mathbb{Z}^+\}$,  \\$A_5=\{(x_1,x_2)\, : \, x_1+2x_2\geq 10,   x_1,x_2\in
 \mathbb{Z}^+\}$,\\
 are not Ameso sets since the  points  $(3,1), (12,4)\in
A_4, $ but $(7,2)\notin A_4,$ and also points $(8,1), (2,4)\in A_5,
$ but $(5,2)\notin A_5.$

\textbf{Example 3:} \ok{Let $D^1=[-20,20] $ and let
$f(x)=\frac{1}{4}x^4-x^3+x,$ $x\in D^1 .$ It is easy to see that
$D^1$ is an Ameso set and $f(x)+f(y)+4\geq
f(\lceil\frac{x+y}{2}\rceil)+f(\lfloor\frac{x+y}{2}\rfloor)$ for
every $x,y\in D^1$. Thus $(D^1, f)$ is an Ameso($4$) pair,
i.e., $C=4$. And the minimization problem of $f(x)$ subject to $ x\in D^1$ is an
Ameso($4$) optimization.}

\ok{Below we state some properties which follow from the definition of an Ameso($C$) pair. The second and third properties can be used for Ameso relaxation   of complicated functions where it is difficult to verify directly that they are an Ameso pair.}

% the following is property 1
\begin{Property_}\label{prop_PlusMinus}
If $(D^n,f)$ is an Ameso($C$) pair then the following inequality
holds for all  $\vec{x}$, $\vec{x}+\vec{a},\vec{x}-\vec{a}\in D^n$
$$f(\vec{x}+\vec{a})+f(\vec{x}-\vec{a})+C\geq 2f(\vec{x}) .$$
\end{Property_}
\begin{pf}%{\bf Property~\ref{prop_PlusMinus}.}
The proof follows from the definition 2 of an Ameso($C$) pair $(D^n,f)$  substituting   $\vec{x}+\vec{a}$ for $\vec{x}$, and  $\vec{x}-\vec{a}\in D^n$ for $\vec{y}$ we obtain:
\bean
f(\vec{x}+\vec{a})+f(\vec{x}-\vec{a})+C
\geq f(\lceil\frac{2\vec{x}}{2}\rceil)+f(\lfloor\frac{2\vec{x}}{2}\rfloor)=2f(\vec{x}).
\eean
%We complete the proof.
\end{pf}

Property~\ref{prop_PlusMinus} implies the relationship between the midpoint convex function and the Ameso function. The more details will be discovered in Property~\ref{prop_midpoint}.

% the following is property 2
\begin{Property_}\label{prop_C2BgC1}
If $(D^n,f)$ is an Ameso($C_1$) pair and $C_2\geq C_1$, then $(D^n,f)$
is an Ameso($C_2$) pair.
\end{Property_}
\begin{pf}%{\bf Property~\ref{Prop_C2BgC1}.}
From an Ameso($C_1$) pair $(D^n,f)$, we have that for every $\vec{x}$, $\vec{y}
\in D^n$,
\bean
f(\vec{x})+f(\vec{y})+C_2\geq
f(\vec{x})+f(\vec{y})+C_1\geq
f(\lceil\frac{\vec{x}+\vec{y}}{2}\rceil)+f(\lfloor\frac{\vec{x}+\vec{y}}{2}\rfloor).
\eean
Hence it is an Ameso($C_2$) pair.
\end{pf}

% the following is Property 3
\begin{Property_}\label{prop_additive}
If $(D^n,f)$ is an Ameso$(C)$ pair, and $(D^n,g)$ is an Ameso$(C')$
pair, then $(D^n,af+bg)$, $a,b\geq 0,$ is an Ameso$(aC+bC')$ pair.
\end{Property_}
\begin{pf}%{\bf Property~\ref{prop_additive}.}
From an Ameso($C$) pair $(D^n,f)$ and  an Ameso($C'$) pair $(D^n,g)$, we have that for every $\vec{x}$, $\vec{y}
\in D^n$ and $a,b\geq 0$,
\bean
&&a[f(\vec{x})+g(\vec{x})]+b[f(\vec{y})+g(\vec{y})]+aC+bC'\\
&&= a[f(\vec{x})+f(\vec{y})+C]+b[g(\vec{x})+g(\vec{y})+C']\\
&& \geq a\big\{f(\lceil\frac{\vec{x}+\vec{y}}{2}\rceil)+f(\lfloor\frac{\vec{x}+\vec{y}}{2}\rfloor)\big\}
+b\big\{g(\lceil\frac{\vec{x}+\vec{y}}{2}\rceil)+g(\lfloor\frac{\vec{x}+\vec{y}}{2}\rfloor)\big\}\\
&&=
a\big\{f(\lceil\frac{\vec{x}+\vec{y}}{2}\rceil)+g(\lceil\frac{\vec{x}+\vec{y}}{2}\rceil)\big\}
+b\big\{f(\lfloor\frac{\vec{x}+\vec{y}}{2}\rfloor)+g(\lfloor\frac{\vec{x}+\vec{y}}{2}\rfloor)\big\}.
\eean
Therefore, $(D^n,af+bg)$ is an Ameso$(aC+bC')$ pair for every $a,b\geq 0$.
\end{pf}

Property~\ref{prop_additive} explains the additive character of the Ameso pair. In the next section, we  further  explore the relation between a convex function and the Ameso pair.

\subsection{Relation Between Ameso Optimization and Convexity}

\ok{In this section we prove that if the domain of a convex function is an Ameso set and
if $f$ is a bounded discrete midpoint convex function, then $(\mathbb{Z}^n,f)$ is an Ameso$(0)$ pair.}
 This  result points to how useful  the Ameso optimization framework can be for discrete optimization problems.

 To start recall the definition of a convex function $f$, for every $\lambda\in[0,1]$,
$$
\lambda f(x)+(1-\lambda)f(y)\geq f(\lambda x+(1-\lambda)y), \, \, x,y\in \mathbb{R}^n.
$$

We first establish the relation between one-dimensional convex function and the Ameso pair in the following Property.

\begin{Property_}\label{prop_1convex}
 For $1\leq i\leq n,$ let $f_i(x)$, $x\in \mathbb{R}$, is a convex function such that there exists a lower bound for each $f_i$. Let $g(\vec{y})=\sum_{i=1}^{n}a_if_i(y_i)$,  $a_i\geq 0$, then
$(\mathbb{Z}^n, g)$ is an Ameso$(0)$ pair.
\end{Property_}
\begin{pf} %{\bf Property~\ref{prop_1convex}.}
Clearly $\mathbb{Z}^n$ is an Ameso set and $g$ has a lower bound. For every $i$, $1\leq i\leq n$ and $x_i,y_i\in \mathbb{Z}$, we will show that
\begin{equation}
f_i(x_i)+f_i(y_i)\geq f_i\big(\big\lceil\frac{x_i+y_i}{2}\big\rceil\big)+f_i\big(\big\lfloor\frac{x_i+y_i}{2}\big\rfloor\big).
\label{eqn_prop_1dim_convex_q01}
\end{equation}

For the case $\frac{x_i+y_i}{2}\in \mathbb{Z}$, we have that $\big\lceil\frac{x_i+y_i}{2}\big\rceil=\big\lfloor\frac{x_i+y_i}{2}\big\rfloor=\frac{x_i+y_i}{2}$. Hence, we have
$$
f_i\big(x_i\big)+f_i\big(y_i\big)\geq 2 f_i\big(\frac{x_i+y_i}{2}\big)=f_i\big(\big\lceil\frac{x_i+y_i}{2}\big\rceil\big)+f_i\big(\big\lfloor\frac{x_i+y_i}{2}\big\rfloor\big).
$$
The inequality comes from the convexity of $f_i$.

Now, if $\frac{x_i+y_i}{2}\notin \mathbb{Z}$, then clearly $x_i\neq y_i$. Let $z^a=min\{x_i,y_i\}$ and $z^b=max\{x_i,y_i\}$. Then, the following three statements are true
\begin{eqnarray}
&z^a-1<z^a<z^a+1\leq z^b, \nonumber\\
&z^b-\frac{z^a+1+z^b}{2}=\frac{z^b-1-z^a}{2}=\frac{z^a-1+z^b}{2}-z^a\geq 0, \nonumber \\
&\frac{z^a+1+z^b}{2}=\big\lceil\frac{z^a+z^b}{2}\big\rceil=\frac{z^a+z^b}{2}+0.5=\big\lfloor\frac{z^a+z^b}{2}\big\rfloor+1.\nonumber
\end{eqnarray}
From $z^b-\frac{z^a+1+z^b}{2}=\frac{z^a-1+z^b}{2}-z^a$, $z^b\geq \frac{z^a-1+z^b}{2}$ and from the  of convexity of  $f_i$, we have
\bean
f_i(z^b)-f_i(\frac{z^a+1+z^b}{2})\geq f_i(\frac{z^a-1+z^b}{2})-f_i(z^a).
\eean
Combining with $\frac{z^a+1+z^b}{2}=\big\lceil\frac{z^a+z^b}{2}\big\rceil$ and $\frac{z^a-1+z^b}{2}=\big\lfloor\frac{z^a+z^b}{2}\big\rfloor$, we have
\begin{eqnarray}
&f_i(z^a)+f_i(z^b)\geq f_i(\frac{z^a+1+z^b}{2})+f_i(\frac{z^a-1+z^b}{2})
=f_i\big(\big\lceil\frac{z^a+z^b}{2}\big\rceil\big)+f_i\big(\big\lfloor\frac{z^a+z^b}{2}\big\rfloor\big).\nonumber
\end{eqnarray}
Also, since $z^a=min\{x_i,y_i\}$ and $z^b=max\{x_i,y_i\}$, we have
$$
f_i(x_i)+f_i(y_i)=f_i(z^a)+f_i(z^b)\geq f_i\big(\big\lceil\frac{z^a+z^b}{2}\big\rceil\big)+f_i\big(\big\lfloor\frac{z^a+z^b}{2}\big\rfloor\big) =f_i\big(\big\lceil\frac{x_i+y_i}{2}\big\rceil\big)+f_i\big(\big\lfloor\frac{x_i+y_i}{2}\big\rfloor\big).
$$

Therefore, we have proved that (\ref{eqn_prop_1dim_convex_q01}) holds for every $i$ and every $x_i, y_i\in \mathbb{Z}$. Hence, we have
\begin{eqnarray}
&&g(\vec{x})+g(\vec{y})=\sum_{i=1}^na_if_i(x_i)+\sum_{i=1}^na_if_i(y_i)=\sum_{i=1}^na_i(f_i(x_i)+f_i(y_i)) \nonumber \\
&&\,\,\,\,\,\,\,\,\,\,\,\,\,\,\,\,\,\,\,\,\,\,\,\,\,\,\,\,\,\,\,\,\,\,\,\,\,\, \geq \sum_{i=1}^n a_i(f_i\big(\big\lceil\frac{x_i+y_i}{2}\big\rceil\big) +f_i\big(\big\lfloor\frac{x_i+y_i}{2}\big\rfloor\big))\nonumber \\
&&\,\,\,\,\,\,\,\,\,\,\,\,\,\,\,\,\,\,\,\,\,\,\,\,\,\,\,\,\,\,\,\,\,\,\,\,\,\,
=g\big(\big\lceil\frac{x_i+y_i}{2}\big\rceil\big)+g\big(\big\lfloor\frac{x_i+y_i}{2}\big\rfloor\big).\nonumber
\end{eqnarray}
Thus, $(\mathbb{Z}^n, g)$ is an Ameso$(0)$ pair.
\end{pf}

The simplest example for Property~\ref{prop_1convex} is the one-dimensional scenario, which setting $n=1$ (see Corollary below).

\begin{Corollary_}\label{cor_1-dimensional_convex_Ameso} If $f$ is one-dimensional convex function with lower bound, then $(\mathbb{Z}, f)$ is an 1-dimensional Ameso$(0)$ pair.
\end{Corollary_}

\begin{Property_}\label{prop_1-Ameso_multi-Ameso}For $1\leq i\leq n,$ let $(\mathbb{Z}, f_i)$ is a 1-dimensional Ameso pair such that there exists a lower bound for each $f_i$. Let $g(\vec{y})=\sum_{i=1}^{n}a_if_i(y_i)$,  $a_i\geq 0$, then
$(\mathbb{Z}^n, g)$ is an Ameso$(0)$ pair.
\end{Property_}
\begin{pf}%{\bf Property~\ref{prop_1-Ameso_multi-Ameso}.}
Clearly, $g$ has a lower bound. For every $i$, $1\leq i\leq n$ and $x_i,y_i\in \mathbb{Z}$, from $(\mathbb{Z}, f_i)$ is a 1-dimensional Ameso pair, we have
\bean
f_i(x_i)+f_i(y_i)\geq f_i\big(\big\lceil\frac{x_i+y_i}{2}\big\rceil\big)+f_i\big(\big\lfloor\frac{x_i+y_i}{2}\big\rfloor\big).
\eean
Hence, we have
\begin{eqnarray}
&&g(\vec{x})+g(\vec{y})=\sum_{i=1}^na_if_i(x_i)+\sum_{i=1}^na_if_i(y_i)=\sum_{i=1}^na_i(f_i(x_i)+f_i(y_i)) \nonumber \\
&&\,\,\,\,\,\,\,\,\,\,\,\,\,\,\,\,\,\,\,\,\,\,\,\,\,\,\,\,\,\,\,\,\,\,\,\,\,\, \geq \sum_{i=1}^n a_i(f_i\big(\big\lceil\frac{x_i+y_i}{2}\big\rceil\big) +f_i\big(\big\lfloor\frac{x_i+y_i}{2}\big\rfloor\big))\nonumber \\
&&\,\,\,\,\,\,\,\,\,\,\,\,\,\,\,\,\,\,\,\,\,\,\,\,\,\,\,\,\,\,\,\,\,\,\,\,\,\,
=g\big(\big\lceil\frac{x_i+y_i}{2}\big\rceil\big)+g\big(\big\lfloor\frac{x_i+y_i}{2}\big\rfloor\big).\nonumber
\end{eqnarray}
Thus, $(\mathbb{Z}^n, g)$ is an Ameso$(0)$ pair.
\end{pf}

From Properties~\ref{prop_additive}, ~\ref{prop_1convex} and ~\ref{prop_1-Ameso_multi-Ameso}, we can compare the convex function and the Ameso as follows.

\begin{Remark_}\label{rm_degenerate}
The Ameso(0) pair not only can be constructed by the 1-dimensional convex function (Property~\ref{prop_1convex}), but also can be obtained by the 1-dimensional Ameso(0) pair (Property~\ref{prop_1-Ameso_multi-Ameso}). Similar to the convexity
property  that the sum of convex functions is a convex function, the sum of Ameso(0) pairs is an Ameso(0) pair (Property~\ref{prop_additive}).
\end{Remark_}

There is one extension of convex function, called discrete midpoint convex function. It can be connected with the Ameso pair also. Recalling the definition of a discrete midpoint convex function as follows,
$$
f(\vec{x})+f(\vec{y})\geq f(\big\lceil\frac{\vec{x}+\vec{y}}{2}\big\rceil)+f(\big\lfloor\frac{\vec{x}+\vec{y}}{2}\big\rfloor), \, \, \, \forall x, y\in \mathbb{Z}^n,
$$
we obtain the following property.

\begin{Property_}\label{prop_midpoint}
If $f$ is a bounded discrete midpoint convex function, then $(\mathbb{Z}^n,f)$ is an Ameso$(0)$ pair.
\end{Property_}
\begin{pf}%{\bf Property~\ref{prop_midpoint}.}
Clearly $\mathbb{Z}^n$ is an Ameso set. Under the assumption $f$ has a lower bound, we also for every $\forall x, y\in \mathbb{Z}^n$,
\bean
f(\vec{x})+f(\vec{y})\geq f(\big\lceil\frac{\vec{x}+\vec{y}}{2}\big\rceil)+f(\big\lfloor\frac{\vec{x}+\vec{y}}{2}\big\rfloor).
\eean
According to Definition~\ref{def_Amesopair}, we have $(\mathbb{Z}^n,f)$ is an Ameso(0) pair. The proof is complete.
\end{pf}

The discrete midpoint convex function is considered as an extension of the convex function on the integer lattice. The Ameso pair can be considered as further extension of the discrete midpoint convex function on a special defined lattice. Furthermore, the condition of midpoint convexity is also extended from $0$ to a non-negative number $C$.

\subsection{Properties of the one-dimensional Ameso($C$) Pair}

\ok{In this section we state and prove   properties of the one dimensional Ameso optimization. The optimization algorithm we propose in this paper is a decomposition algorithm i.e., it is based on  finding  solutions to simpler one dimensional problems, see also    \cite{GGW2011}, \cite{KV1987}. We start with the following lemma which shows us the form of an Ameso set. According to this lemma an Ameso set can be expressed only in the form $[a,b]$.}

\begin{Lemma_}\label{lem_1-dimensional_integer}
 An one-dimensional set $M\subset \mathbb{Z}$ is an Ameso set if and only if it can be expressed
as $[x_s,x_t]$,  with $x_s< x_t$.
\end{Lemma_}

\begin{pf} First consider a  set $M=[x_s,x_t] $.
Then for all $x,y \in M$, we have that $x,y \in [x_s,x_t]$ and $x,y$
are integers.  Also, it is easy to see that
\bean
min\{x,y\}\leq\left\lfloor\frac{x+y}{2}\right\rfloor\leq\left\lceil\frac{x+y}{2}\right\rceil\leq
max\{x,y\}.
\eean
 Therefore, $
\left\lceil\frac{x+y}{2}\right\rceil,\left\lfloor\frac{x+y}{2}\right\rfloor\in
M$. That is, $M$ is an one-dimensional Ameso set.

Conversely, we will prove that every one-dimensional Ameso set,  $D^1$, can be
expressed as $[x_s,x_t] $ where $x_s < x_t$ are
integers.

First by relabeling we can write  $D^1 =
\{a_1,a_2,...,a_m\}$ where $a_{i-1}<a_i$, $a_i\in Z$ for all $i$.

From $a_{i-1}<a_i$, $a_i\in Z$ for all $i$, we have $a_{i-1}\leq a_i-1$.

We next show that $a_{i}-a_{i-1}=1$ for every $i$. We assume there exists $\hat{i}$ with $a_{\hat{i}}-a_{\hat{i}-1}\geq 2$.  Hence,
\bean
\lfloor\frac{a_{\hat{i}}+a_{\hat{i}-1}}{2}\rfloor
\geq
\lfloor\frac{2a_{\hat{i}-1}+2}{2}\rfloor
=\lfloor a_{\hat{i}-1}+1 \rfloor=a_{\hat{i}-1}+1;\\
\mbox{ and } \lfloor\frac{a_{\hat{i}}+a_{\hat{i}-1}}{2}\rfloor
\leq
\lfloor\frac{2a_{\hat{i}}-2}{2}\rfloor
=\lfloor
a_{\hat{i}}-1\rfloor=a_{\hat{i}}-1.
\eean
Therefore,
\bean
a_{\hat{i}-1}<\lfloor\frac{a_{\hat{i}}+a_{\hat{i}-1}}{2}\rfloor<a_{\hat{i}}.
\eean
Combining with $D^1 =
\{a_1,a_2,...,a_m\}$ where $a_i < a_{i+1}$ for all $i$, we have
$\lfloor\frac{a_{\hat{i}}+a_{\hat{i}-1}}{2}\rfloor \notin D^1$ by
its construction. It contradicts the definition of the Ameso set.

Thus, for all $i$, $a_i-a_{i-1}=1$ then $D^1$  can be expressed by
$[a_1,a_m] $. The  proof is complete.
\end{pf}

To avoid trivial cases in the sequel, we assume that the Ameso
sets under study are not empty or singletons. Further  according
to Lemma 1, $D^1=[x_s,x_t] $ where $x_s$ is the smallest possible value of the independent variable and $x_t$  is the biggest value of the independent variable. To avoid the trivial case, we assume $x_s<x_t$.  In the following discussion in one dimension, we will directly use $[x_s,x_t]$.

Next we discuss properties of the one-dimensional Ameso($C$) pair $([x_s,x_t],f)$. In the following discussion, we first narrow the local optimizer into a small range (Lemma 2) and then explore how to extend the local optimizer to the global optimizer by providing certain conditions (Lemma 3). Furthermore, we provide relaxed conditions of local optimizer to be a global optimizer (Theorem 1).

% the following is lemma 2
\begin{Lemma_}\label{lem_rightnarrow_local}
\ok{If there exist}
 $ x^0\in [x^s,x^t]$ and $b\in\mathbb{Z}^+$ such that
$[x^0,x^0+b]\subseteq [x^s,x^t]$ and the  following conditions hold
\begin{itemize}
  \item[(a)] $f(x^0)= min_{y\in [x^0,x^0+b]}f(y)$;
  \item[(b)] $f(x^0)+C\leq max_{y\in [x^0,x^0+b] }f(y)$.
\end{itemize}
then $\arg^*\max_{y\in [x^0,x^0+b]}f(y)\in
(x^0+\frac{b}{2},x^0+b] $ where
$\arg^*$ means the largest optimizer if multiple optimizers exist.
\end{Lemma_}
\begin{pf}Let $z=\arg^*\max_{y\in [x^0,x^0+b] }f(y)$. We need to prove that $z\in
(x^0+\frac{b}{2},x^0+b] $ under the two given conditions in the statement. From the condition (b) of the statement and definition of $z$  we have that $f(x^0)+C\leq f(z)$.

If $C=0$, from Property~\ref{prop_PlusMinus} we have that for all $x, f(x+1)+f(x-1)\geq 2f(x)$, which can be written as
\bean
f(x+1)-f(x)\geq f(x)-f(x-1), \, \, \forall x.
\eean
Combining with the condition $f(x^0)=
min_{y\in [x^0,x^0+b] }f(y)$, we have $f$ is increasing in this interval $x\geq x^0$.

Therefore, $x_0+b=\arg^*\max_{y\in [x^0,x^0+b] }f(y)$ for the case
$C=0$.

Now we assume that $C > 0$. We need to show that $z\in (x^0+\frac{b}{2},x^0+b]$. We assume
$z\in[x_0, x^0+\frac{b}{2}]$ and then we obtain a contradiction.

First, from $f(x^0)+C\leq f(z)$ and $C>0$, we have
$z\neq x^0$.

Second, from the assumption we have $z \in
(x^0,x^0+\frac{b}{2}] $. Hence
\bean
2z-x^0\in
(x^0,x^0+b]
\eean
\bea\label{eqn_lem_rightnarrow_local_q1}
\mbox{ and \, \, } 2z-x^0=z+z-x^0>z+0=z \Rightarrow 2z-x^0>z.
\eea
From $z=\arg^*\max_{y\in [x^0,x^0+b] }f(y)$, we also have
\bea\label{eqn_lem_rightnarrow_local_q2}
f(2z-x^0)<f(z).
\eea
Hence,
\bea\label{eqn_lem_rightnarrow_local_q3}
f(x^0)+f(z)+C>f(x^0)+f(2z-x^0)+C\geq 2f(z).
\eea
The first inequality comes from Eq.(\ref{eqn_lem_rightnarrow_local_q2}). The second inequality comes from the Ameso($C$) pair $(D^1,f)$.

Eq.(\ref{eqn_lem_rightnarrow_local_q3}) can be written as $f(x^0)+C> f(z)$, which contradicts
with  the condition (a) in the statement $f(x^0)+C\leq f(z)$.

Hence, $z\in (x^0+\frac{b}{2},x^0+b]$ for the case $C>0$.

We complete the proof.
\end{pf}

Lemma~\ref{lem_rightnarrow_local} shows that the location of the local optimal point in a domain can be narrowed into a half length of domain if we find that the discrete function is an Ameso pair and satisfies the two listed conditions. We will apply the result to the right side.

% the following is lemma 3
\begin{Lemma_}\label{lem_narrow_right}
\ok{If there exist} $ x^0\in [x_s,x_t]$ and $b\in\mathbb{Z}^+$ such that
$[x^0,x^0+b]\subseteq [x_s,x_t]$ and the  following conditions hold
\begin{align*}
f(x^0)&= min_{y\in [x^0,x^0+b] }f(y),\\
f(x^0)+C &\leq max_{y\in [x^0,x^0+b] }f(y) ,
\end{align*}
then $ f(x^0)= min_{y\in [x^0,x^t] }f(y)$.
\end{Lemma_}
\begin{pf} First, if $x^0+b=x_t$, then the statement follows by assumption.

For $x^0+b<x_t$, we will use mathematical induction to prove it.

From $ [x^0,x^0+b+1]\supset [x^0,x^0+b]$,  we have that
\bea
\max_{y\in [x^0,x^0+b+1] }f(y)\geq \max_{y\in
[x^0,x^0+b] }f(y).\label{eqn_lem_narrow_right_q1}
\eea
From the condition (b) of the statement $ f(x^0)+C\leq max_{y\in
[x^0,x^0+b] }f(y)$, it follows that $ f(x^0)+C\leq max_{y\in
[x^0,x^0+b+1] }f(y)$.

Next we will show that $f(x^0)= min_{y\in [x^0,x^0+b+1 ]}f(y)$.

Let $z=\arg^*\max_{y\in[x^0,x^0+b]}f(y)$. By Lemma \ref{lem_rightnarrow_local} with the condition $f(x^0)+C\leq
max_{y\in [x^0,x^0+b]}f(y)$, we have that $z\in (x^0+\frac{b}{2},x^0+b]$.

Hence,
\bean
2z\in
(2x^0+b,2x^0+2b]
\, \, \, \Rightarrow \, \, \,
2z-(x^0+b+1)\in (x^0-1,x^0+b-1].
\eean
Since $ 2z-(x^0+b+1)$ is integer,
we have that
\bean
2z-(x^0+b+1)>x^0-1
\, \, \, \Rightarrow \, \,\,
2z-(x^0+b+1)\geq x^0.
\eean
Then,
\bean
2z-(x^0+b+1)\in [x^0,x^0+b-1].
\eean Therefore,
\bea\label{eqn_lem_narrow_right_q2}
f(2z-(x^0+b+1))\leq f(z).
\eea
Hence
\bea
f(x^0+b+1)+f(z)+C
\geq f(x^0+b+1)+f(2z-(x^0+b+1))+C
\geq 2f(z). \label{eqn_lem_narrow_right_q3}
\eea
The first inequality comes from Eq.(\ref{eqn_lem_narrow_right_q2}). The second inequality comes from the Ameso($C$) pair $([x_s,x_t],f)$.
Thus,
\bean
f(x^0+b+1)\geq f(z)-C \geq f(x^0)
\eean
The first inequality comes from Eq.(\ref{eqn_lem_narrow_right_q3}). The second inequality comes from the condition of the statement $f(x^0)+C \leq max_{y\in [x^0,x^0+b] }f(y)=f(z)$.
Hence, $f(x^0)\leq
f(x^0+b+1)$.

Also, from
$[x^0,x^0+b+1] =([x^0,x^0+b])\cup\{x^0+b+1\}$,
\bean
\min_{[x^0,x^0+b+1]}f(y)=\min_{y\in {[x^0,x^0+b]} }f(y)=f(x^0).
\eean

From above argument we show that for the interval $[x^0,x^0+b+1]$, both $f(x^0)=
min_{y\in {[x^0,x^0+b]+1} }f(y)$ and $f(x^0)+C\leq max_{y\in
{[x^0,x^0+b]+1} }f(y)$ hold.

Follow a similar argument, we can obtain $[x^0,x^0+b+k+1]$ from the condition $f(x^0)= min_{y\in [x^0,x^0+b+k] }f(y)$
and $f(x^0)+C\leq max_{y\in [x^0,x^0+b+k]}f(y)$.

Therefore, from mathematical induction, we have
$$f(x^0)= min_{y\in [x^0,x_t] }f(y)$$ and $$f(x^0)+C\leq
max_{y\in [x^0,x_t] }f(y),$$ hold. This completes the
proof.
\end{pf}

Lemma~\ref{lem_narrow_right} narrows the search of the left side of an interval for a local optimizer. We now extend the result to the global optimizer.

% the following is theorem 1
\begin{Theorem_}\label{th_GlobalOptimizer_right}
\ok{If there exist} $ x^0\in [x_s,x_t]$ and $b\in\mathbb{Z}^+$ such that
$[x^0,x^0+b]\subseteq [x_s,x_t]$ and the  following conditions hold
\begin{align*}
f(x^0) &= min_{y\in [x_s, x^0+b] }f(y), \\
f(x^0)+C&\leq max_{y\in [x^0, x^0+b] }f(y),
\end{align*}
then $ f(x^0)= min_{y\in [x_s,x_t]}f(y)$.
\end{Theorem_}
\begin{pf}
From $ f(x^0)= \min_{y\in ([x_s,x^0+b])}f(y)$ and $x^0\in [x_s,x^0], x^0\in [x^0,x^0+b]$, we have that $
f(x^0)= min_{y\in [x_s,x^0] }f(y)$, and $f(x^0)= min_{y\in
{[x^0,x^0+b]} }f(y)$.

Now, from the condition $ f(x^0)+C\leq max_{y\in
{[x^0,x^0+b]}}f(y)$ and according to Lemma~\ref{lem_narrow_right}, we have that $ f(x^0)= min_{y\in {[x^0,x_t]} }f(y)$. Also, $f(x^0)= min_{y\in
[x_s,x^0] }f(y)$. Therefore $ f(x^0)= min_{y\in [x_s,x_t]}f(y)$.
\end{pf}

Theorem~\ref{th_GlobalOptimizer_right} shows that if there exists an interval $[a, b]\subseteq[x_s,x_t]$ in the domain where $f(a)$ is the minimum of $f$ in the interval $[x_s, b]$ and follows certain condition, then this local minimum is the global minimum, i.e. $f(a)$ is the minimum in $[x_s, x_t]$.

\ok{Below we give a corollary which provides a way to find an interval in the domain with the properties we derived  above. }

% the following is corollary 1
\begin{Corollary_}\label{col_z_xprime_rightrange}
\ok{If there exist} $ x'\in [x_s,x_t], z\in (x',x_t] $ with
$f(z)-f(x')\geq C$, then $$min_{y\in
[x_s,z] }f(y)=min_{y\in [x_s,x_t]}f(y).$$
\end{Corollary_}
\begin{pf}
Let $f(x^0)= min_{y\in [x_s,z] }f(y)$, we
need to prove that $min_{y\in [x_s,x_t]}f(y)=f(x^0)$.

Let $f(\hat{x})=min_{y\in [x',z] }f(y)$ and $x'\leq \hat{x}\leq z$, then
 $f(\hat{x})=min_{y \in
[\hat{x},z] }f(y)$ and $f(\hat{x})+C\leq f(x')+C\leq
f(z)\leq max_{y \in [\hat{x},z] }f(y)$.

According to Lemma~\ref{lem_narrow_right},
\bea\label{eqn_col_z_xprime_rightrange_q1}
f(\hat{x})=min_{y\in [\hat{x},x_t] }f(y)
\eea

From $ [x_s,z]\supseteq [\hat{x},z]$ and $f(x^0)=min_{y\in
[x_s,z] }f(y)$, $f(\hat{x})=min_{y\in
[\hat{x},z] }f(y)$, we have that $ f(x^0)\leq f(\hat{x})$. Using (\ref{eqn_col_z_xprime_rightrange_q1}) and $\hat{x}\leq z$, we have that $f(x^0)\leq min_{y\in
[\hat{x},x_t] }f(y)$ and
\bea
f(x^0)=min_{y\in
[x_s,z] }f(y)\leq min_{y\in
[x_s,\hat{x}] }f(y).\label{eqn_col_z_xprime_rightrange_q2}
\eea
Using Eq.(\ref{eqn_col_z_xprime_rightrange_q1}) and Eq.(\ref{eqn_col_z_xprime_rightrange_q2}),  we have that $
f(x^0)=min_{y\in [x_s,x_t]}f(y)$. 
That is, \\$min_{y\in
[x_s,z] }f(y)=min_{y\in [x_s,x_t]}f(y)$.
\end{pf}

\ok{In the sequel we state the analogous  lemmas, theorem and corollary for the case where there exists an interval $[a,b]$ in the domain where now $b$ is the minimum of $f$. Therefore, all the proofs are omitted since they are similar analogous to previous ones. }

% the following is lemma 4
\begin{Lemma_}\label{lem_leftnarrow_local}
\ok{If there exist} $ x^0\in [x_s,x_t]$ and $b\in\mathbb{Z}^+$ such that
$[x^0-b,x^0]\subseteq {[x_s,x_t]}$ with the  following conditions holding
\begin{align*}
f(x^0) &= min_{y\in [x^0-b,x^0] }f(y), \\
f(x^0)+C&\leq max_{y\in [x^0-b,x^0] }f(y),
\end{align*}
then $$min\{w : f(w)=max_{y\in [x^0-b,x^0]}f(y)\} \in
[x^0-b,x^0-\frac{b}{2}) .$$
\end{Lemma_}

% the following is lemma 5
\begin{Lemma_}\label{lem_narrow_left}
\ok{If there exist} $ x^0\in [x_s,x_t]$ and $b\in\mathbb{Z}^+$ such that
$[x^0-b,x^0]\subseteq {[x_s,x_t]}$ with the  following conditions holding
\begin{align*}
f(x^0) &= min_{y\in [x^0-b,x^0] }f(y) \\
f(x^0)+C &\leq max_{y\in [x^0-b,x^0] }f(y)
\end{align*}
then $ f(x^0)= min_{y\in [x_s,x^0] }f(y)$.
\end{Lemma_}

% the following is theorem 2
\begin{Theorem_}\label{th_GlobalOptimizer_left}
\ok{If there exist} $ x^0\in [x_s,x_t]$ and $b\in\mathbb{Z}^+$ such that
$[x^0-b,x^0]\subseteq {[x_s,x_t]}$ with the  following conditions holding
\begin{align*}
f(x^0) &= min_{y\in [x^0-b,x_t] }f(y) \\
f(x^0)+C&\leq max_{y\in [x^0-b,x^0] }f(y)
\end{align*}
then $ f(x^0)= min_{y\in [x_s,x_t]}f(y)$.
\end{Theorem_}

% the following is corollary 2
\begin{Corollary_}\label{col_z_xprimer_leftrange}
\ok{If there exist} $ x'\in [x_s,x_t], z\in[x_s,x') $ with
$f(z)-f(x')\geq C$, then $$min_{y\in [z,x_t] }f(y)= min_{y\in [x_s,x_t]}f(y).$$
\end{Corollary_}

\ok{The following theorem is important since it proves that if we can find an interval in which there is a local minimum under some conditions then this is global minimum. To do that we use the properties of the intervals we defined in the above discussion.}

% the following is theorem 3
\begin{Theorem_}
\ok{If there exist} $ x^0\in [x_s,x_t]$ and $ b_1,b_2\in\mathbb{Z}^+,$ such
that $ [x^0-b_2,x^0+b_1]\subseteq {[x_s,x_t]}$ with the following
conditions:
\begin{align*}
f(x^0)&= min_{y\in [x^0-b_2,x^0+b_1] }f(y) ,\\
f(x^0)+C &\leq max_{y\in [x^0,x^0+b_1] }f(y), \\
f(x^0)+C&\leq max_{y\in [x^0-b_2,x^0] }f(y),
\end{align*}
then $ f(x^0)= min_{y\in [x_s,x_t]}f(y)$.
\end{Theorem_}
\begin{pf}From $f(x^0)= min_{y\in [x^0-b_2,x^0+b_1] }f(y)$, we have
\bean
f(x^0)= min_{y\in [x^0-b_2,x^0] }f(y)= min_{y\in [x^0,x^0+b_1] }f(y)
\eean

Also, from $ f(x^0)= min_{y\in
{[x^0,x^0+b_1]} }f(y),$ $f(x^0)+C\leq max_{y\in
{[x^0,x^0+b_1]} }f(y)$ and Lemma~\ref{lem_narrow_right}, we have $ f(x^0)= min_{y\in
{[x^0,x_t]} }f(y)$.

Finally, from $ f(x^0)=
min_{y\in [x^0-b_2,x^0] }f(y),$ $f(x^0)+C\leq max_{y\in
[x^0-b_2,x^0]}f(y)$ and Lemma~\ref{lem_narrow_left}, we have $ f(x^0)= min_{y\in
[x_s,x^0] }f(y)$.

Hence,
\bean
f(x^0)=\min_{y\in
[x^0,x_t] }f(y)=\min_{y\in
[x_s,x^0]\bigcup[x^0,x_t] }f(y)
=\min_{y\in
[x_s,x_t]}f(y).
\eean
The proof is complete.
\end{pf}

\ok{Now, we state a corollary which shows how to find the interval with the above property.}

% the following is corollary 3
\begin{Corollary_}\label{col_z_xprimer_range}
\ok{If there exist} $x', z_s, z_t\in [x_s,x_t],$ such that  $ z_s<x',$ $z_t>x'$
with $f(z_s)-f(x')\geq C$ and $f(z_t)-f(x')\geq C$

then $$min_{y\in [z_s,z_t] }f(y)=min_{y\in [x_s,x_t]}f(y).$$
\end{Corollary_}
\begin{pf} Let $f(x^0)=min_{y\in [z_s,z_t] }f(y)$,
we need to prove $min_{y\in [x_s,x_t]}f(y)=f(x^0)$.

Since from
$f(x^0)=min_{y\in [z_s,z_t] }f(y)=min_{y\in
[x^0,z_t] }f(y)$ and the condition in the statement, we have $ f(x^0)+C\leq f(x')+C\leq
f(z_t)\leq max_{y\in [x^0,z_t] }f(y)$. Thus, $
f(x^0)=min_{y\in [x^0,x_t] }f(y)$ according Lemma~\ref{lem_narrow_right}.

Also, from $f(x^0)=min_{y\in [z_s,z_t] }f(y)=min_{y\in
[z_s,x^0] }f(y)$ and the condition in the statement, we have $f(x^0)+C\leq f(x')+C\leq
f(z_s)\leq max_{y\in [z_s,x^0] }f(y)$.
Thus,
$f(x^0)=min_{y\in [x_s,x^0] }f(y)$ according Lemma~\ref{lem_narrow_left}.

Finally,
\bean
f(x^0)=min_{y\in [x^0,x_t] }f(y)=min_{y\in [x_s,x^0] }f(y)= min_{y\in [x_s,x_t]}f(y).
\eean
We complete the proof.
\end{pf}

Through the above discussion, we have narrowed the computation of
the optimal solution of any one-dimensional Ameso($C$) optimization problem
into the identification of  some intervals given in the above
theorems.

\ok{The following Example implements the above corollary to minimize a function with domain an Ameso set. }

\textbf{Example 4:} $f(x)=\frac{1}{4}x^4-x^3+x$, $x\in
[-20,20] $ is an Ameso(4) set from Example 2. At the
same time, we have $f(1)-f(3)=4\geq 4, $ $f(4)-f(3)=7.75\geq 4$
and $f(3)=min\{f(i): i=1,2,3,4\}, $ then $f(3)=min\{f(i): i\in
[-20,20] \}.$

\ok{As we mentioned above corollary 4 is useful to obtain an algorithm which can solve an one dimension Ameso optimization problem. Therefore, we have the following algorithm.}

    \ok{\begin{center}
    \small
    \fbox{
    \begin{tabular}{ll}
    \multicolumn{2}{c}{\textbf{one-dimensional Ameso($C$) optimization algorithm}}\\
    \textbf{Input}: & Ameso($C$) optimization problem: \\
    &$minimize$ $
    f(x);$ $subject \ to$ $ x\in [x_s,x_t]  $\\
    \textbf{Step 1}: & $A=\phi, l^+=0, l^-=0,$ select a point  $l^0\in [x_s,x_t]$,$l^*=l^0, S=l^0$,\\
     &  then calculate and define as
    $f(l^0)=\min  \ f(x)$, \\
    \textbf{Step 2}: & Update $l^+=min\{x\, : \, x>l^0, x\in  [x_s,x_t]- A\}$,\\
    &  hence $f(l^+)=\min \ f(x)$, \\
    \textbf{Step 3}:
    & If $f(l^+)- f(l^*)\geq C$, go to \textbf{step 4};\\
    & If $0 \leq f(l^+)- f(l^*)< C$, then $A=A\cup \{l^+\}$ go to \textbf{step 2}.\\
    & If $ f(l^+)- f(l^*)< 0$, then $A=A\cup \{l^+\}$, $l^*=l^+$ go to \textbf{step 2}.\\
    & If $\{x\, : \, x\in [x_s,x_t]-A, x>l_0\}=\phi$, go to \textbf{step 4};\\
    \textbf{Step 4}:& Update $l^-$ to be  any integer satisfying $f(l^-)=max_{\{l\, : \, l\in A, l<l^*\}}f(l)$;\\
    & If $f(l^-)- f(l^*)\geq C$, go to \textbf{output};\\
    \textbf{Step 5}: & Update $l^-=max\{x\, : \, x<l^0, x\in  [x_s,x_t]- A\}$,\\
    &  then calculate and define as $f(l^-)=\min \ f(x)$ \\
    \textbf{Step 6}:
    & If $f(l^-)- f(l^*)\geq C$, go to \textbf{output};\\
    & If $0 \leq f(l^-)- f(l^*)< C$, then $A=A\cup \{l^-\}$ go to \textbf{step 5}.\\
    & If $ f(l^-)- f(l^*)< 0$, then $A=A\cup \{l^-\}$, $l^*=l^-$ go to \textbf{step 5}.\\
    & If $\{x\, : \, x\in [x_s,x_t]-A, x<l_0\}=\phi$, go to \textbf{output};\\
    \textbf{Output}: & $A$, $l^*$, $f(l^*)$\\
    \end{tabular}
    }
    \end{center}}

\ok{Below we give an example to show how the algorithm can be implemented.}
%\pagebreak[4]

\textbf{Example 5:} Define a function $f(x)$ where the range of
$x$ is the set $$\{1, 2, 3, 4, 5, 6, 7, 8, 9, 10, 11, 12, 13, 14,
15, 16, 17, 18, 19, 20, 21, 22, 23, 24, 25, 26, 27, 28, 29, 30,
31\}$$ with corresponding values
$$f(x)=7,9,7,8,7,8,9,8,7,8,9,8,7,8,6,7,4,5,6,7,7,6,7,8,9,10,11,9,8,7,8
,$$ i.e., $f(1) =7,$ $f(2) =9,\ldots,$ $f(31) =8.$ The graph of
$f$ is given figure 1. The problem of minimizing $f$ over its
range is an one-dimensional Ameso(7) problem. \\

%\begin{center}
%\epsfig{file=1.eps,height=8cm}\\
%Figure 1
%\end{center}

\begin{figure}[!h]
%\vspace*{-2cm}
\begin{center}
\includegraphics[width=0.4\textwidth]{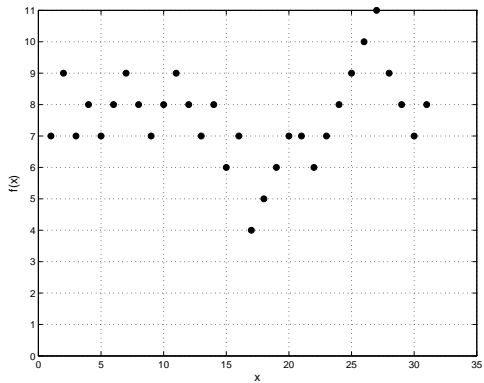}
\end{center}
%\vspace*{+2cm}
\caption{Function $f(x)$ of Example 5.}\label{Fig_1}
\end{figure}

If  $l^0=13$, then we can find the global minimum
point $l^*=17$ after doing computations  involving only the set of
points in  the set $\{1,\cdots,27\}$. Indeed the algorithm will
work as follows: first it will  search to  the right of $13$ and
it will stop at point $l^+=27$  because $f(27)-f(17)=7.$ The
minimum value in the interval $\{13,\cdots,27\}$ is $f(17)=4.$

Then the algorithm will compare the maximum value $f(14)=8$ of
$f$ in the interval $\{13,\cdots,16\}$ and since $f(14)-f(17)=4
<7$, it will continue with the following step.

Now, the algorithm will search  to the left of $13$ and it will stop at
point $l^-=1$  because there is no $x\in \{1,\cdots,13\}$
satisfying $f(x)-f(17)\geq 7.$  The algorithm has computed
the global minimum point $l^*=17,$ without considering points in
the set $\{28,\cdots,31\}.$

\ok{From the definition of an Ameso($C$) optimization problem one can see that  any discrete optimization
problem over a finite set can be transformed into an Ameso($C$)
optimization problem when $C\geq 2*(max\{f(x)\} - min\{f(x)\}),$
 and its domain is an Ameso set. However, such a large $C$ is meaningless,
because in order to satisfy the conditions of stopping in the above algorithm we need to check the whole domain if $C> max\{f(x)\} - min\{f(x)\}$. Therefore, it is obvious that an Ameso($C$)  optimization can be preferred over other discrete optimization methods only if $C \leq max\{f(x)\} - min\{f(x)\}$. In such cases the above algorithm can narrow down the computations dramatically and as we showed in Example 5.
For instance, in the Example 5 the function $f$ is an Ameso(C) problem, for any
$C\geq 7.$ If we use the algorithm with  $C=8$, and starting point at $l^0=13$, we have going though its steps as above that it will search the whole interval $\{1,\cdots, 31\}$. This happened since in that case $C=8>7=max\{f(x)\} - min\{f(x)\}$. Another question that raises from the above example is which $C$ to choose since we know that if a problem is Ameso$(C)$ it is also Ameso$(aC)$, $a>0$. As we showed in Example 5, the number of computations depends on the starting point and $C$.}

Notice that for an Ameso($C$) problem, $C$ is a critical value to the problem. In the algorithm, from the condition $f(l^+)-f(l^-)\geq C$, we can see that $C$ defines the search range. Generally, the algorithm takes more time to search the optimal solution for higher values of $C$.

\subsection{Properties of multi-dimensional Ameso($C$) Pair}

In this section we introduce the multi-dimensional Ameso optimization problem and we discuss its properties. Using these properties we show an algorithm for the multi-dimensional case which is based on the decomposition analysis and the one-dimensional Ameso optimization problem.

Due to that it is very complex problem, we add the following requirement to the domain of the optimization.

\begin{Assumption_}\label{ass_multiD}
 For the multi-dimensional problem,  assume that the domain can be written as follows.
\bea \label{eqn_ass_multiD}
D^n=\{(x_1,\cdots, x_i): x_i\in D_i\}
\eea
where $D_i$ is the domain of the $i-th$ dimension variable   is an Ameso set.
\end{Assumption_}

For a multi-dimensional Ameso$(C)$ optimization problem:
$minimize$ $ f(\vec{x});$ $subject \ to $ $\vec{x}\in D^n$,  we construct another problem: fixing the value of $j$ dimensions, denotes as $i_1, \cdots, i_j$, where $i_k\in [1,n]$ for any $1\leq k\leq j$ satisfying $k\neq l, i_k\neq i_l$. For
fixed $i_1,...,i_j (n\geq j)$  we start with the following
definition.\\

% define the condition domain
\begin{Definition_}

\begin{enumerate}
  \item {\bf The domain} $\Delta^j_{i_1,...,i_j}$
 is the set:
$$\Delta^j_{i_1,...,i_j}=\{(x_{i_1},\cdots,x_{i_j}):\exists
(x_1',...,x_n')\in D^n \mbox{ with } x_{i_{k}}'=x_{i_k} \   ,\ k=1,...,j\}.$$

which is the domain of  $x_{i_1}$, $\cdots$ $x_{i,j}$ in the original domain $D^n.$

Under Assumption~\ref{ass_multiD},  
$\Delta^j_{i_1,...,i_j}=D_{i_1}\times D_{i_2}\times\cdots\times D_{i_j}$.

\item {\bf The conditional domain:} For    fixed
$x_{i_{1}}^0,...,x_{i_{j}}^0$ we call the  set:
$$\Gamma_{x_{i_{1}}^0,...,x_{i_{j}}^0}^{n-j}=\{(x_1,...,x_n)\in
D^n: \ x_{i_{k}}=x_{i_{k}}^0 \   ,\ k=1,...,j\}$$
the conditional domain of $D^n$ given $x_{i_{1}}^0,...,x_{i_{j}}^0.$ 
\item {\bf The conditional   function}: $f_{i_1,...,i_j}^*: \Delta^j_{i_1,...,i_j} \longrightarrow \Re$

$$f_{i_1,...,i_j}^*(x_{i_{1}},...,x_{i_{j}})=min_{ \vec{y}\in
 \Gamma_{x_{i_{1}},\ldots,x_{i_{j}}}^{n-j}} f(\vec{y}).$$
 \item {\bf  The conditional pair} of Ameso(C) pair  to be
  the pair: $(\Delta_{i_1,...,i_j}^j,f_{i_1,...,i_j}^*).$

\end{enumerate}
\end{Definition_}

We first establish the following   essential property.

% the following is property
\begin{Property_}\label{pro_conditionalPair}
Under Assumption~\ref{ass_multiD}, the conditional pair $(\Delta_{i_1,...,i_j}^j,f_{i_1,...,i_j}^*)$
of an n-dimensional Ameso($C$) pair $(D^n,f)$ is a j-dimensional
Ameso($C$) pair.
\end{Property_}
\begin{pf} At first we prove $\Delta_{i_1,...,i_j}^j$ is  an
Ameso set. And $\forall (x_{i_1},x_{i_2},...,x_{i_j}),
(y_{i_1},y_{i_2},...,y_{i_j}) \in \Delta_{i_1,...,i_j}^j$, we have for all $\vec{x}'=(x_1',...,x_n')\in D^n,
x_{i_k}'=x_{i_k},\forall k=1,...j$ and a vector
$\vec{y}'=(y_1',...,y_n')\in D^n, y_{i_k}'=y_{i_k},\forall
k=1,...j$.

From $D^n$ is an Ameso set, we have $\vec{x}',\vec{y}'\in
D^n\Rightarrow\left\lceil\frac{\vec{x}'+\vec{y}'}{2}\right\rceil,\left\lfloor\frac{\vec{x}'+\vec{y}'}{2}\right\rfloor\in
D^n$.

Thus, there exists an
$\left\lceil\frac{\vec{x}'+\vec{y}'}{2}\right\rceil\in D^n $ with
$\left\lceil\frac{x_{i_k}'+y_{i_k}'}{2}\right\rceil=\left\lceil\frac{x_{i_k}+y_{i_k}}{2}\right\rceil,
\forall k=1,...j$ and

 $\left\lfloor\frac{\vec{x}'+\vec{y}'}{2}\right\rfloor\in D^n$
 with
$\left\lfloor\frac{x_{i_k}'+y_{i_k}'}{2}\right\rfloor=\left\lfloor\frac{x_{i_k}+y_{i_k}}{2}\right\rfloor,
\forall k=1,...j$.

From  $\Delta_{i_1,...,i_j}^j=D_{i_1}\times D_{i_2}\cdots D_{i_2}$, we have
\begin{equation}\label{eq1}
(\lceil\frac{\vec{x}_{i_1}+\vec{y}_{i_1}}{2}\rceil,...,\lceil\frac{\vec{x}_{i_j}+\vec{y}_{i_j}}{2}\rceil),
(\lfloor\frac{\vec{x}_{i_1}+\vec{y}_{i_1}}{2}\rfloor,...,\lfloor\frac{\vec{x}_{i_j}+\vec{y}_{i_j}}{2}\rfloor)
\in \Delta_{i_1,...,i_j}^j
\end{equation}
So $\Delta_{i_1,...,i_j}^j$ is an Ameso set.

From the fact that function $f$ has a lower bound and the $f(\vec{x}_0)=f_{i_1,...,i_j}^*(x_{i_1},x_{i_2},...,x_{i_j})$ is in the range of $f$, the value $f(\vec{x}_0)=f_{i_1,...,i_j}^*(x_{i_1},x_{i_2},...,x_{i_j})$ must be bigger than the lower bound of $f$, for every $(x_{i_1},x_{i_2},...,x_{i_j})$. Thus, $f(\vec{x}_0)=f_{i_1,...,i_j}^*(x_{i_1},x_{i_2},...,x_{i_j})$ has a lower bound for each $(x_{i_1},x_{i_2},...,x_{i_j})$.

Next we  show the third condition of the definition. For any
 $ (x_{i_1},x_{i_2},...,x_{i_j}),
(y_{i_1},y_{i_2},...,y_{i_j}) \in \Delta_{i_1,...,i_j}^j$,
$\exists \vec{x}_0\in D^n,\vec{y}_0\in D^n$,
$f(\vec{x}_0)=f_{i_1,...,i_j}^*(x_{i_1},x_{i_2},...,x_{i_j}),f(\vec{y}_0)=f_{i_1,...,i_j}^*(y_{i_1},y_{i_2},...,y_{i_j}).$
Since  $(D^n,f)$ is an Ameso($C$) pair,  we have
\begin{center}
$\begin{array}{rl}
&f_{i_1,...,i_j}^*(x_{i_1},x_{i_2},...,x_{i_j})+f_{i_1,...,i_j}^*(y_{i_1},y_{i_2},...,y_{i_j})+C\\
&=f(\vec{x}_0)+f(\vec{y}_0)+C\geq
f(\lceil\frac{\vec{x}_0+\vec{y}_0}{2}\rceil)+
f(\lfloor\frac{\vec{x}_0+\vec{y}_0}{2}\rfloor)\\
&\geq
f_{i_1,...,i_j}^*(\lceil\frac{\vec{x}_{i_1}+\vec{y}_{i_1}}{2}\rceil,...,\lceil\frac{\vec{x}_{i_j}+\vec{y}_{i_j}}{2}\rceil)+f_{i_1,...,i_j}^*(\lfloor\frac{\vec{x}_{i_1}+\vec{y}_{i_1}}{2}\rfloor,...,\lfloor\frac{\vec{x}_{i_j}+\vec{y}_{i_j}}{2}\rfloor)\\
&\mbox{(because of (\ref{eq1}) and the definition of
$f_{i_1,...,i_j}^*(\cdot)$)}.
\end{array}$
\end{center}

Now, we have a $j-dimensional$ Ameso set $\Delta_{i_1,...,i_j}^j$ and a function with lower bound satisfies the condition of Ameso pair $f_{i_1,...,i_j}^*$. Therefore, $(\Delta_{i_1,...,i_j}^j,f_{i_1,...,i_j}^*)$ is
also an Ameso($C$) pair, and $f_{i_1,...,i_j}^*$ is a $j$-dimensional
function for every $(x_{i_1},\cdots, x_{i_j})\in \Delta_{i_1,...,i_j}^j$.
\end{pf}

Now, using the above property we can construct an algorithm which obtains the minimum of a multi-dimensional Ameso optimization problem. Therefore, we establish the following theorem.

% the following is theorem 4
\begin{Theorem_}

The solution of any Ameso optimization problem can be found using
the Ameso Recursive Procedure(ARP), described in the table below.

\end{Theorem_}

\begin{center}
\small
\fbox{
\begin{tabular}{ll}
\multicolumn{2}{c}{\textbf{Ameso Recursive
Procedure (ARP)}}\\
\textbf{Input}: & Ameso($C$) Optimization problem: \\
&$minimize$ $
f(\vec{x});$ $subject \ to$ $ \vec{x}\in D^n  $\\
\textbf{Step 1}: & $A=\phi, l^+=0, l^-=0,$ select a point  $l^0\in\Delta^1_{n}$,$l^*=l^0, S=l^0$,\\
 &  then calculate
$f_{n}^*(l^0)=min_{(x_1,\ldots,x_{n-1},l^0)\in
\Gamma_{l^0}^{n-1} }f(x)$  \ \dag , \\
\textbf{Step 2}: & Update $l^+=min\{x\, :  \, x>l^0, x\in  \Delta^1_{n}- A\}$,\\
&  then calculate $f_{n}^*(l^+)=min_{(x_1,\ldots,x_{n-1},l^+)\in
\Gamma_{l^+}^{n-1} }f(x)$  \dag\\
\textbf{Step 3}:
& If $f_{n}^*(l^+)- f_{n}^*(l^*)\geq C$, go to \textbf{step 4};\\
& If $0 \leq f_{n}^*(l^+)- f_{n}^*(l^*)< C$, then $A=A\cup \{l^+\}$ go to \textbf{step 2}.\\
& If $ f_{n}^*(l^+)- f_{n}^*(l^*)< 0$, then $A=A\cup \{l^+\}$, $l^*=l^+$ go to \textbf{step 2}.\\
& If $\{x\, :  \, x\in\Delta^1_{n}-A, x>l_0\}=\phi$, go to \textbf{step 4};\\
\textbf{Step 4}:& Update $l^-$ to be  any integer satisfying $f_{n}^*(l^-)=max_{\{l\, :  \, l\in A, l<l^*\}}f_{n}^*(l)$;\\
& If $f_{n}^*(l^-)- f_{n}^*(l^*)\geq C$, go to \textbf{output};\\
\textbf{Step 5}: & Update $l^-=max\{x\, :  \, x<l^0, x\in  \Delta^1_{n}- A\}$,\\
&  then calculate $f_{n}^*(l^-)=min_{(x_1,\ldots,x_{n-1},l^-)\in
\Gamma_{l^-}^{n-1} }f(x)$ \dag\\
\textbf{Step 6}:
& If $f_{n}^*(l^-)- f_{n}^*(l^*)\geq C$, go to \textbf{output};\\
& If $0 \leq f_{n}^*(l^-)- f_{n}^*(l^*)< C$, then $A=A\cup \{l^-\}$ go to \textbf{step 5}.\\
& If $ f_{n}^*(l^-)- f_{n}^*(l^*)< 0$, then $A=A\cup \{l^-\}$, $l^*=l^-$ go to \textbf{step 5}.\\
& If $\{x\, :  \, x\in\Delta^1_{n}-A, x<l_0\}=\phi$, go to \textbf{output};\\
\textbf{Output}: & $A$, $l^*$, $f_{n}^*(l^*)$\\
 %&
%$(x_1*,\ldots,x_{n-1}*,l^*)$:
% $f_{n}^*(l^*)=f(x_1*,\ldots,x_{n-1}*,l^*)$\\
\multicolumn{2}{l}{\dag Note  : the computation of $f_{n}^*(l^0),$
$f_{n}^*(l^+)$ and $f_{n}^*(l^-)$ above requires }\\
\multicolumn{2}{l}{ solving $(n-1)$-dimensional Ameso($C$)
optimization problems.}\\
\end{tabular}
}
\end{center}

\begin{pf}  
The proof is easy to complete using  Property~\ref{pro_conditionalPair} and
Corollary~\ref{col_z_xprimer_range}.
\end{pf}

Below, we give an example with a 2-dimensional  Ameso optimization problem and we implement the {\bf ARP} algorithm we provided above.

\textbf{Example 6:} Consider the function
$f(x_1,x_2)=88e^{\frac{1}{x_1}}+99
e^{\frac{2}{x_2}}+\, :  \, \frac{sin(x_1x_2)}{2}\, :  \, , $  and the domain
$D^2=\{( x_1,x_2) \, :  \,  \ x_i =1,2,\cdots,100, \  i=1,2 \}$. The
problem of minimizing $f(x)$ over $D^2$  is a two-dimensional Ameso($1$)
optimization.

%Without loss of generality start with

 Assume we pick  $x_2=80$ as the starting
point.  Then $l^0=80$. The {\bf ARP} will calculate
$f_{2}^*(80)=min_{(x_1,80)\in \Gamma_{80}^{1}}f(x_1,80)=190.4220, $
 where $\Gamma_{80}^{1}=\{(x_1,x_2)\, :  \, x_1=1,\cdots,100;  x_2=80\}.$ Then it will begin to search in $x_2$'s increasing side,
 i.e, $l^+=81,82,83,...$.
It will keep on calculating $f_{2}^*(l^+)=min_{(x_1,l^+)\in
\Gamma_{l^+}^{1}}f(x_1,l^+). $  In the end it will stop at
$l^+=100,$ and the the minimum point in the interval $x_2 \in \{80,
\ldots,100\}$ is $f_{2}^*(97)=min_{(x_1,97)\in
\Gamma_{97}^{1}}f(x_1,97)=f(97, 97)=190.0093. $ %i.e. here
%corresponding $x_1=97$ with $f_{2}^*(97)=f(97, 97).$

Next we calculate the difference of the maximum point  in the set
$\{80, \cdots, 96\}$ and $f_{2}^*(97)$, that is: $
f_{2}^*(81)-f_{2}^*(97)<1=C$. Hence the ARP  will next search in the
direction of the decreasing side of $x_2 (x_2< 80),$  it will
compute $f_{2}^*(l^-)$ for  $l^-=79,78,\ldots $ and  it stop when
$l^-=66$ with $f_{2}^*(66)-f_{2}^*(97)=191.1640-190.0093\geq 1.$

Now  we can say that $x_1^*=97, x_2^*=97$ is a global minimum with
value $f(x_1^*,x_2^*)=190.0093.$

We illustrate the function $f$ and the {\bf ARP} technique with the
following figures. In figure 2 we graph the function for all
points in its domain $D^2.$ In figure 3 we graph the function for
a smaller part of its domain. In figure 4 we plot the function
only for points that were involved  in our {\bf ARP} technique.
Note that the number of points in figure 4 is significantly
smaller than the points in Figure 3 and figure 2. In figure 5 we
plot the function $f_{2}^*(x_2), $ that is calculated in our {\bf
ARP} technique. It is  an one-dimensional Ameso(1) function.

\begin{center}
\epsfig{file=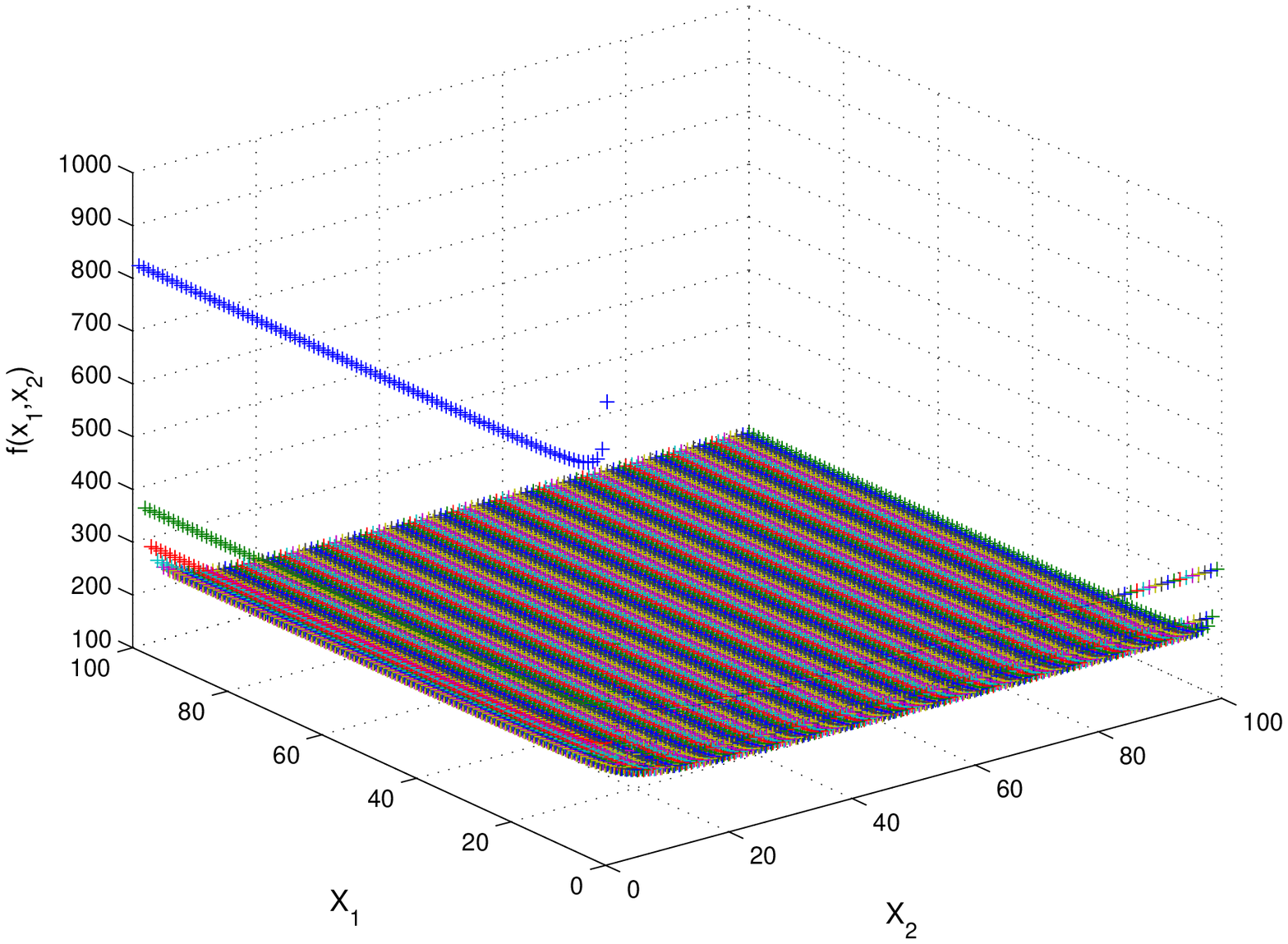,height=6cm}\\
Figure 2
\end{center}

\begin{center}
\epsfig{file=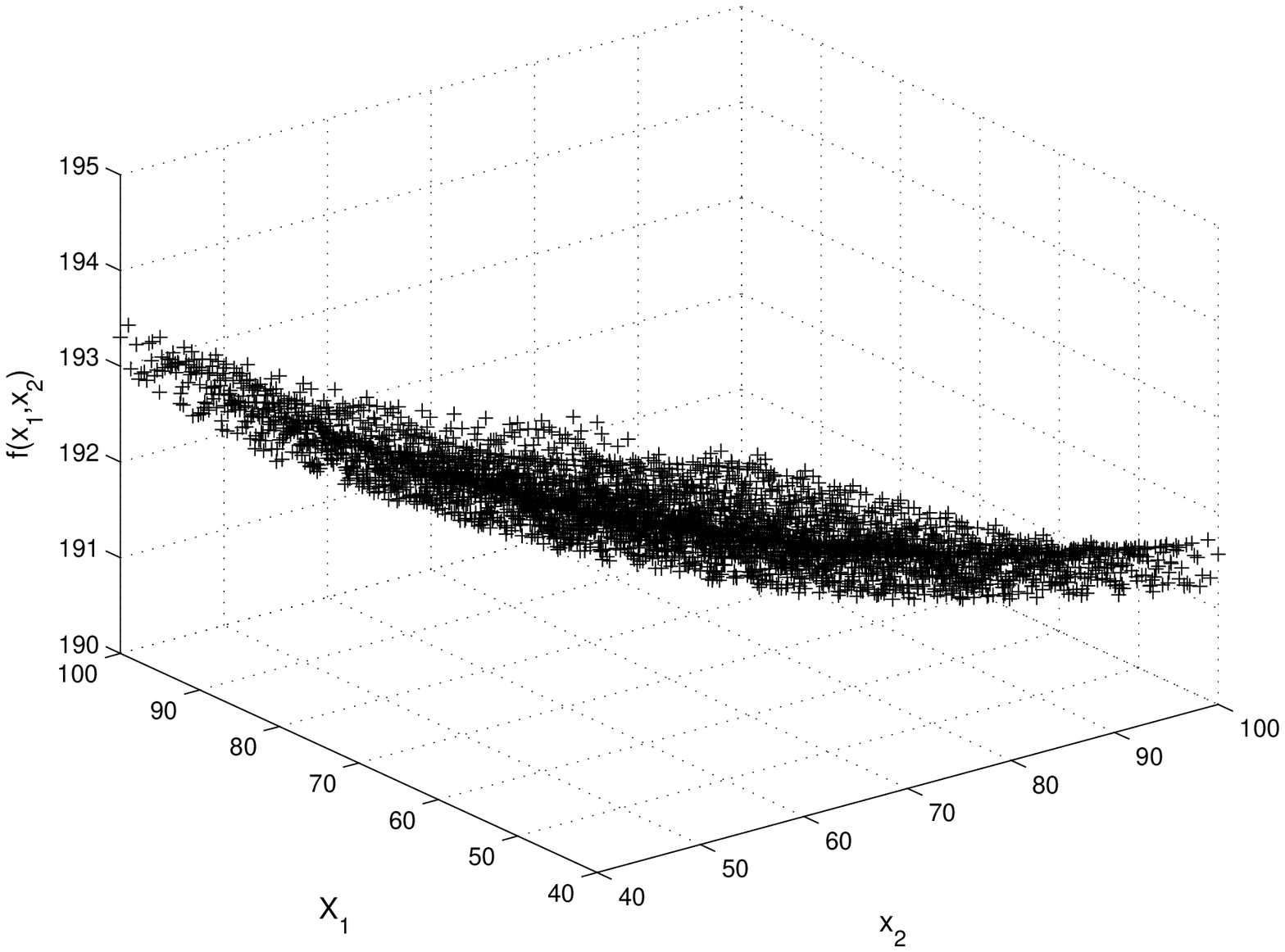,height=6cm}\\
Figure 3
\end{center}

\begin{center}
\epsfig{file=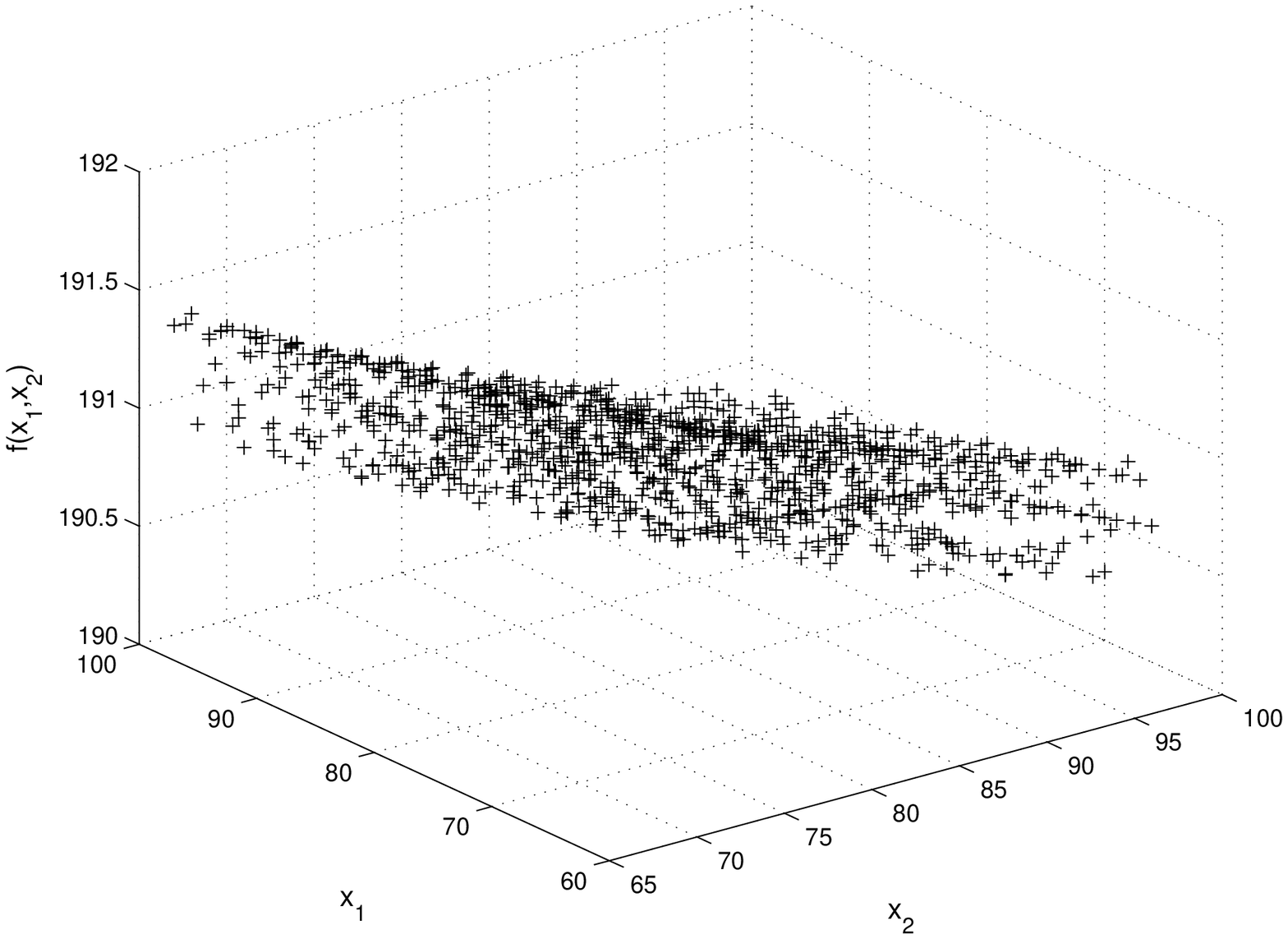,height=6cm}\\
Figure 4
\end{center}

\begin{center}
\epsfig{file=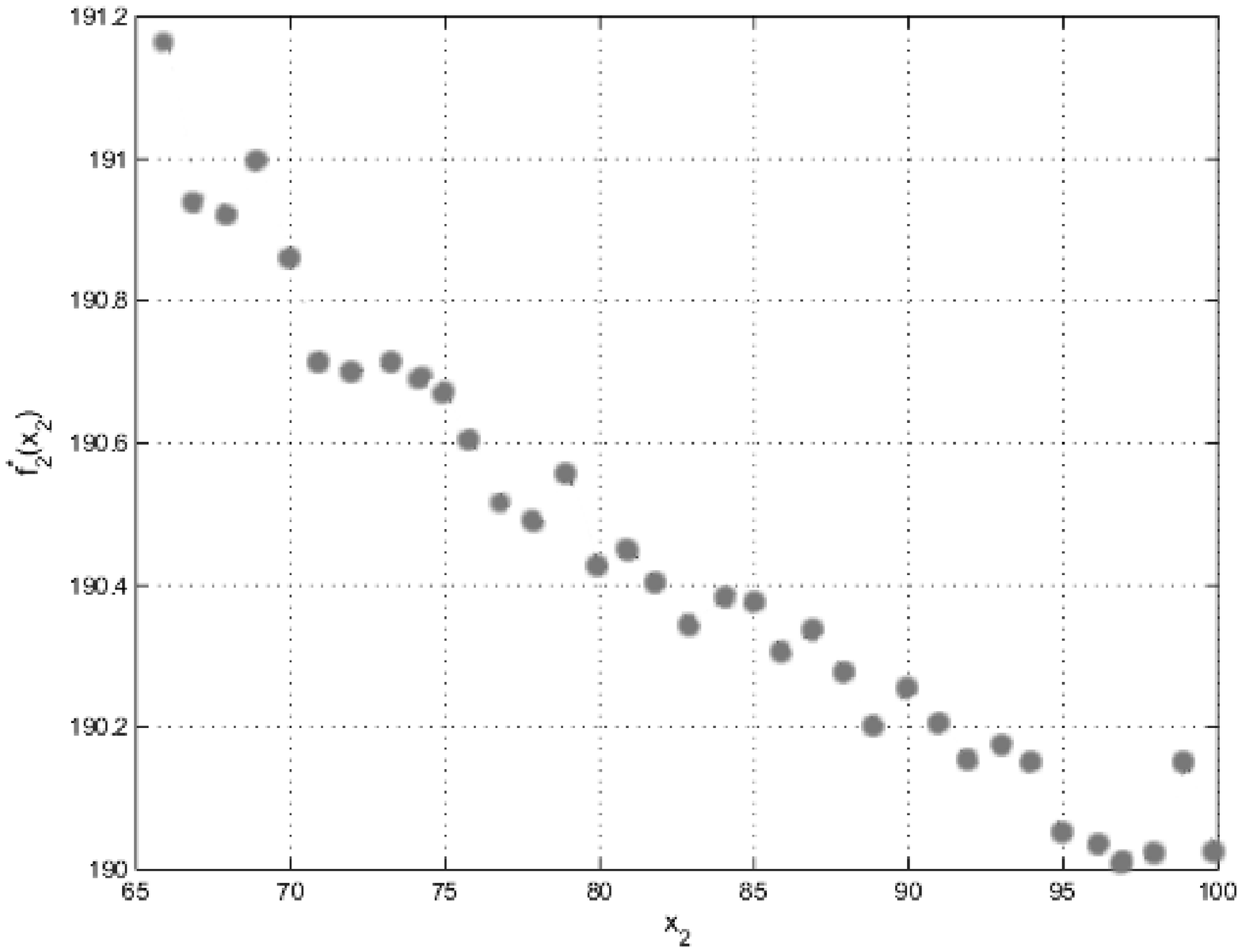,height=6cm}\\
Figure 5
\end{center}

\begin{Remark_}
In Example 6, if we choose  any point $x_2 $
greater than  $65 $ as starting point, then the {\bf ARP} will
stop after a  search of all $l^+\in \{x_2+1,\cdots, 100\}$ and
$l^-\in \{x_2-1, \cdots, 66\}$. If choose as starting point $x_2 $
one that is smaller than $66 ,$  then the {\bf ARP} will stop
after a  search of all $l^+\in \{x_2+1,\cdots, 100\}$, and it will
not do a  search in the left side of $x_2$.
\end{Remark_}
\begin{Remark_} From Examples 5 and 6, we can see that the
choice of  starting points and the form of the functions
themselves determine  the complexity of {\bf ARP}. For  an
arbitrary integral optimization problem, if one  can establish
that  it is an Ameso($C$) optimization for a suitable number  $C,$
then one can use the  {\bf ARP} to find the optimal  solution
without necessarily searching all points in the domain.
\end{Remark_}
%\subsection{Ameso(1) Optimization and Midpoint Convexity}

\section{An Application of Ameso(C)}
In this section we consider a well known version of the
  knapsack problem and show  it is an Ameso optimization problem. 
A  manufacture has $W$  units of an item for shipping to its retailer. There are $N=3$ shipping options as follows. Option $i$: each package with capacity $w_i \in \mathbf{Z}^+$ and cost $c_i \in \mathbf{Z}^+$; $i=1,2,3$. 
It is assumed that: 
$w_1<w_2<w_3$, $c_1<c_2<c_3$, and $\frac{c_1}{w_1}>\frac{c_2}{w_2}>\frac{c_3}{w_3}.$

Let $z_i$  denote the number of packages using option $i$, where $z_i\in \mathbf{Z}^+\bigcup\{0\},$ $i=1,2$. It is  assumed that  any 
 remaining capacity:  $\lceil\frac{W-w_1z_1-w_2z_2}{w_3}\rceil$  will be assigned to  Option 3. 
  
   From $\frac{c_1}{w_1}>\frac{c_2}{w_2}>\frac{c_3}{w_3}$, we have that Option 1 has the highest unit shipping cost and Option 3 is the cheapest one. The manufacturer has an external constraint of 
   upper limit for the number of packages   $z_i$, such as $z_iw_i\leq \frac{W}{2}, \, i=1,2$.

Hence, the total cost of manufacturer can be written as

\bea\label{eqn_app_obj}
C(z_1,z_2)=c_1z_1+c_2z_2+c_3\lceil\frac{W-w_1z_1-w_2z_2}{w_3}\rceil, \, \, \, z_iw_i\leq \frac{W}{2}, \, i=1,2.
\eea

The target is to minimize the total cost under the constraint. Therefore, the ideal structure of the objective function $C(z_1,z_2)$ in (\ref{eqn_app_obj}) is to be convex function or a function related to convex function in discrete version. The first discretizing convex function is the midpoint convex function which we have discussed (Property~\ref{prop_midpoint}). One can prove that $C$ is not midpoint convex function. Therefore, we can prove that it is an Ameso optimization problem. To do that we establish the following Lemma which is necessary for proving the Ameso optimization problem cf. Theorem \ref{th_app_ameso}.

\begin{Lemma_}\label{lem_app_integer}
\begin{itemize}
  \item[(i)] If $n_1,n_2\in \mathbf{Z}$, then
$n_1+n_2=\lfloor\frac{n_1+n_2}{2}\rfloor+\lceil\frac{n_1+n_2}{2}\rceil$.
  \item[(ii)] If $x_1,x_2\in \mathbf{R}$, then
$\lceil x_1\rceil+\lceil x_2\rceil \leq
\lceil x_1+x_2\rceil +1$ and
$$\lfloor x_1\rfloor+\lfloor x_2\rfloor \leq \lfloor\frac{x_1+x_2}{2}\rfloor+\lceil\frac{x_1+x_2}{2}\rceil
  \leq \lceil x_1\rceil+\lceil x_2\rceil
.$$
%\item[(iii)] If $n\in \mathbf{z}$, then
%$
%\lceil\frac{n}{2}\rceil+\lfloor\frac{n}{2}\rfloor=n$.
\end{itemize}
\end{Lemma_}
\begin{pf}%{\bf lem_app_integer.}
\emph{(i)} 
If both $n_1$ and $n_2$ are even numbers (or both are odd numbers), then
\bean
\lfloor\frac{x_1+x_2}{2}\rfloor+\lceil\frac{x_1+x_2}{2}\rceil
=\frac{x_1+x_2}{2}+\frac{x_1+x_2}{2}=x_1+x_2.
\eean

If $n_1$ is odd number and $n_2$ is even, then
\bean
\lfloor\frac{x_1+x_2}{2}\rfloor+\lceil\frac{x_1+x_2}{2}\rceil
=\frac{x_1+x_2}{2}-0.5+\frac{x_1+x_2}{2}+0.5=x_1+x_2.
\eean

\emph{(ii)} For $x_1, x_2 \in\mathbf{R}$, let
\bean
z_1^-:=\lfloor x_1\rfloor, \,\, a_1:=x_1-z_1, \, \, z_1^+:=\lceil x_1\rceil, \,\, b_1:=z_1^+-x_1;\\
z_2^-:=\lfloor x_2\rfloor, \,\, a_2:=x_2-z_2, \, \, z_2^+:=\lceil x_2\rceil, \,\, b_2:=z_2^+-x_2.
\eean
By definition, $a_1, b_1, a_2, b_2\geq 0$. Thus,
\bean
\lceil x_1\rceil+\lceil x_2\rceil
&&=z_1^+ +z_2^+
=\lceil z_1^+ + z_2^+ -1\rceil+1\\
&&
\leq \lceil z_1^+ + z_2^+ -1 +(1-b_1-b_2)\rceil +1
=\lceil x_1+x_2\rceil +1.
\eean
The inequality comes from $-1<1-b_1-b_2$. Also,
\bean
\lfloor x_1\rfloor+\lfloor x_2\rfloor
&&=z_1^- + z_2^-
=\lfloor\frac{z_1^-+z_2^-}{2}\rfloor+\lceil\frac{z_1^-+z_2^-}{2}\rceil \\
&& \leq \lfloor\frac{z_1^-+z_2^-+a_1+a_2}{2}\rfloor+\lceil\frac{z_1^-+z_2^-+a_1+a_2}{2}\rceil\\
&&= \lfloor\frac{x_1+x_2}{2}\rfloor+\lceil\frac{x_1+x_2}{2}\rceil\\
&& \leq \lfloor\frac{x_1+x_2+b_1+b_2}{2}\rfloor+\lceil\frac{x_1+x_2+b_1+b_2}{2}\rceil\\
&&=\lfloor\frac{z_1^+ +z_2^+}{2}\rfloor+\lceil\frac{z_1^++z_2^+}{2}\rceil
=z_1^+ +z_2^+=\lceil x_1\rceil+\lceil x_2\rceil.
\eean
The second and fifth equality come from the result of \emph{(i)}. Also, the first inequality comes from $a_1, a_2\geq 0$,\footnote{when $a_1>0$ or $a_2>0$, the strictly inequality holds.} 
%The third equality comes from the definition of $x_1=z_1^-+a_1, \, x_2=z_2^-+a_2$. 
and the second inequality comes from $b_1, b_2\geq 0$.\footnote{when $b_1>0$ or $b_2>0$, the strictly inequality holds.} %The fourth equality comes from definition of $x_1=z_2^- - b_1, \, x_2=z_2^+ - b_2$.
%\emph{To see Part iii,} from the result of part i, we have
%\bean
%\lfloor\frac{n}{2}\rfloor+\lceil \frac{n}{2}\rceil
%=\lfloor\frac{n-1+1}{2}\rfloor+\lceil \frac{n-1+1}{2}\rceil
%=n-1+1=n.
%\eean
We complete the proof.
\end{pf}

%In the view of the proof of the following Lemma~\label{lem_app_integer}, we have $\pi$ is not mid-point convex function. The next discretizing convex function is our Ameso(C) function referred in the paper. In the next Theorem, we will prove it is an Ameso function.

\begin{Theorem_}\label{th_app_ameso}
Minimization of $C(z_1,z_2)$ subject to $z_i\in \mathbf{Z}^+\bigcup\{0\},$ $i=1,2$ is an Ameso($c_3$) optimization problem.
\end{Theorem_}
\begin{pf}%{\bf Theorem~\ref{th_app_ameso}.}
When $z_iw_i\leq \frac{W}{2}, \, i=1,2$,  we have that the domain is an Ameso set.

Next, we prove that $C$ has the necessary property to be Ameso($c_3$) pair. If a $(\{\mathbf{Z}^+\bigcup\{0\}\}^2,C)$ is Ameso($c_3$) pair then minimization of $C$ is an Ameso($c_3$) optimization problem. For every $(z_1, z_2), (z_1^\prime, z_2^\prime)\in \{\mathbf{Z}^+\bigcup\{0\}\}^2$, we have that
\bean
&&\,\, C(z_1,z_2)+ C(z_1^{\prime},z_2^{\prime})+c_3\\
&&=c_1(z_1+z_1^{\prime})+c_2(z_2+z_2^{\prime})
+c_3\big(\lceil\frac{W-w_1z_1-w_2z_2}{w_3}\rceil+\lceil\frac{W-w_1z_1^{\prime}-w_2z_2^{\prime}}{w_3}\rceil\big)+c_3\\
&&\geq c_1\big(\lfloor\frac{z_1+z_1^{\prime}}{2}\rfloor+\lceil\frac{z_1+z_1^{\prime}}{2}\rceil\big)
+c_2\big(\lfloor\frac{z_2+z_2^{\prime}}{2}\rfloor+\lceil\frac{z_2+z_2^{\prime}}{2}\rceil\big)
+c_3\big\lceil \frac{W-w_1z_1-w_2z_2}{w_3}+\frac{W-w_1z_1^\prime-w_2z_2^\prime}{w_3}
\big\rceil \\
&&\,\,\,\,\,\,+c_3\\
&&=c_1\big(\lfloor\frac{z_1+z_1^{\prime}}{2}\rfloor+\lceil\frac{z_1+z_1^{\prime}}{2}\rceil\big)
+c_2\big(\lfloor\frac{z_2+z_2^{\prime}}{2}\rfloor+\lceil\frac{z_2+z_2^{\prime}}{2}\rceil\big)
+c_3\big\lceil \frac{2W-w_1(z_1+z_1^{\prime})-w_2(z_2+z_2^\prime)}{w_3}
\big\rceil +c_3\\
&&= c_1\big(\lfloor\frac{z_1+z_1^{\prime}}{2}\rfloor+\lceil\frac{z_1+z_1^{\prime}}{2}\rceil\big)
+c_2\big(\lfloor\frac{z_2+z_2^{\prime}}{2}\rfloor+\lceil\frac{z_2+z_2^{\prime}}{2}\rceil\big)
\\
&&\,\,\,\,\,\,\,\,\,\,\,\,\,\,\,\,\,\,\,\,\,\,\,\,\,\,\,\,\,\,\,\,\,\,\,\,\,\,\,\,\,\,\,\,\,\,
+c_3\bigg(\big\lceil \frac{2W-w_1\big(
\lfloor\frac{z_1+z_1^{\prime}}{2}\rfloor+
\lceil\frac{z_1+z_1^{\prime}}{2}\rceil\big)
-w_2\big(
\lfloor\frac{z_2+z_2^{\prime}}{2}\rfloor+
\lceil\frac{z_2+z_2^{\prime}}{2}\rceil\big)}{w_3}
\big\rceil+1\bigg)\\
&&=c_1\big(\lfloor\frac{z_1+z_1^{\prime}}{2}\rfloor+\lceil\frac{z_1+z_1^{\prime}}{2}\rceil\big)
+c_2\big(\lfloor\frac{z_2+z_2^{\prime}}{2}\rfloor+\lceil\frac{z_2+z_2^{\prime}}{2}\rceil\big)
\\
&&\,\,\,\,\,\,\,\,\,\,\,\,\,\,\,\,\,\,\,\,\,\,\,\,\,\,\,\,\,\,\,\,\,\,\,\,\,\,\,\,\,\,\,\,\,\,
+c_3\bigg(\big\lceil \frac{W-w_1\lfloor\frac{z_1+z_1^{\prime}}{2}\rfloor
-w_2\lfloor\frac{z_2+z_2^{\prime}}{2}\rfloor}{w_3}
+\frac{W-w_1\lceil\frac{z_1+z_1^{\prime}}{2}\rceil
-w_2\lceil\frac{z_2+z_2^{\prime}}{2}\rceil}{w_3}
\big\rceil+1\bigg)
\eean
\bean
&&\geq
c_1\big(\lfloor\frac{z_1+z_1^{\prime}}{2}\rfloor+\lceil\frac{z_1+z_1^{\prime}}{2}\rceil\big)
+c_2\big(\lfloor\frac{z_2+z_2^{\prime}}{2}\rfloor+\lceil\frac{z_2+z_2^{\prime}}{2}\rceil\big)
\\
&&\,\,\,\,\,\,\,\,\,\,\,\,\,\,\,\,\,\,\,\,\,\,\,\,\,\,\,\,\,\,\,\,\,\,\,\,\,\,\,\,\,\,\,\,\,\,
+c_3\bigg[\big\lceil \frac{W-w_1\lfloor\frac{z_1+z_1^{\prime}}{2}\rfloor
-w_2\lfloor\frac{z_2+z_2^{\prime}}{2}\rfloor}{w_3}\big\rceil
+\big\lceil\frac{W-w_1\lceil\frac{z_1+z_1^{\prime}}{2}\rceil
-w_2\lceil\frac{z_2+z_2^{\prime}}{2}\rceil}{w_3}
\big\rceil\bigg]
\\
&&=C(\lfloor\frac{z_1+z_1^{\prime}}{2}\rfloor,\lfloor\frac{z_2+z_2^{\prime}}{2}\rfloor)
+C(\lceil\frac{z_1+z_1^{\prime}}{2}\rceil,\lceil\frac{z_2+z_2^{\prime}}{2}\rceil).
\eean
The first inequality comes from the fact that $
\lceil x_1+x_2\rceil\leq \lceil x_1\rceil+\lceil x_2\rceil$ for every $x_1, x_2\in \mathbf{R}$. %The second equality comes from the simplify the notation. 
Also, the third equality and the second inequality come from Lemma~\ref{lem_app_integer}. Thus,   the proof is complete.
\end{pf}

We have established  that the problem is an Ameso($c_3$) optimization problem.  This work 
can be easily generalized to the case of $k$ possible shipping options. % 

%*******
%At the same time, we also notice that the problem is a revised knapsack problem. As we all know knapsack problem is a NP-complete problem. We also know that many NP-hard problem can be transformed to knapsack NP-complete problem.
%
%The relation between P and NP is still a challenge problem in the discrete computation theory. The definition of Ameso problem provide other point of view to approach this type of problem. We are no longer stick on the description of the problem itself. We try to look at the functions and looking for the answers in the integer domain.
%******

 \section{Conclusions}  

In this paper we introduced the class of Ameso($C$) optimization
problems and  established the following:

\begin{enumerate}
  \item For one-dimensional Ameso optimization problems there are
simple  to verify  optimality conditions at any optimal point cf. 
\textbf{Theorems 1,2,3 and corresponding corollaries}.

  \item The conditional pair of Ameso($C$) pair is still an Ameso($C$)
pair cf. \textbf{Property 5}.

  \item We have
constructed  the {\bf ARP} algorithm that solves multi-dimensional Ameso optimization
problems without necessarily performing complete enumeration cf. 
\textbf{ ARP}.

\end{enumerate}

We showed that  one can establish  that a problem is Ameso($C$) optimization 
in the same way one  proves convexity or midpoint convexity. Further, one can see that an Ameso(C) optimization is a relaxed convex model, since   convex cases correspond to an  
 Ameso(0) optimization problem.

As we mentioned in Section 2.2, in a general setting it is not possible 
 establish bounds for the constant $C$. However, in specific problems  such as in Section 3, 
  a   problem specific bound for  the constant $C$ can be obtained  by direct  model analysis. 
  %In some cases  one can try to prove that the problem is  is an Ameso(1) optimization problem.
%  
Further,     for models with    functions where it is not possible  to prove convexity or midpoint convexity it may still be possible  to show that the  more relaxed property of the Ameso(1) optimization holds. To do that is easier than to prove convexity and at the same time there is the {\bf ARP}  algorithm that can compute  the minimum point without requiring complete enumeration.

{\bf Acknowledgments}  We acknowledge support for this work from the National Science Foundation, NSF grant CMMI-16-62629.

\bibliographystyle{elsarticle-num}

\begin{thebibliography}{las}

%
\bibitem{AMO}
Ahuja, R.K., Magnanti, T.L., Orlin, J.B. (1993). ``Network Flows—Theory, Algorithms and Applications",  \emph{Prentice-Hall}, Englewood Cliffs.
%
\bibitem{BDP}
Bertsekas, D. P. and Tsitsiklis, J. N. (2003). ``Parallel and
Distributed Computation: Numerical Methods'', \emph{LIDS Technical
Reports}, MIT, Boston Ma.
%
\bibitem{BG}
Bertsimas, D.J. and van Ryzin,  G.J.  (1993). ``Stochastic and
Dynamic Vehicle Routing with General Interarrival and Service Time
Distributions", \emph{Advances in Applied Probability}, 25, pp.
947-978.
%
\bibitem{BGJ}
Bertsimas, D.J. and  van Ryzin, G.J. (1990). ``An Asymptotic
Determination of the Minimum Spanning Tree and Minimum Matching
Constants in Geometrical Probability," \emph{Operations Research
Letters}, 9, pp. 223-231.
%
\bibitem{DEL}
De Loera, J.A., Hemmecke, R., Koppe, M. (2013). ``Algebraic and Geometric Ideas in the Theory of Discrete Optimization", \emph{SIAM}, Philadelphia.
%
\bibitem{Denardo1982}
Denardo, E. V. (1982). ``Dynamic Programming: Models and Applications", \emph{NJ, Englewood Cliffs:Prentice-Hall}. 

\bibitem{DEV}
Denardo, E. V.,  Huberman, G. and  Rothblum, U. G. (1982). ``Optimal
Locations on a Line Are Interleaved'',  \emph{Operations Research},
30, pp. 745-759.
%
\bibitem{HKW}
Hemmecke, R., Koppe, M., Lee, J., Weismantel, R. (2010). ``Nonlinear integer programming",  \emph{In: Junger, M., et al. (eds.) 50 Years of Integer Programming 1958-2008: From The Early Years and State-of-the-Art}, Chapter 15, pp. 561-618, Springer, Berlin.
%
\bibitem{GilmoreGomony1961}
Gilmore, P.C. ,and Gomony, R. E. (1961). ``A linear programming approach to the cutting stock problem'', \emph{Operations Research}, 9, pp. 849-859.

\bibitem{GilmoreGomony1963}
Gilmore, P.C. ,and Gomony, R. E. (1963). ``A linear programming approach to the cutting
stock problem. Part II.'', \emph{Operations Research}, 11, pp. 863-888.

\bibitem{GilmoreGomony1966}
Gilmore, P.C. ,and Gomony, R. E.(1966). ``The theory and computation of knapsack functions'', \emph{Operations Research}, 14, pp. 1045-1074.

\bibitem{GGW2011}
Gittins, J., Glazebrook, K., and  Weber, R. (2011). `` \emph{Multi-armed bandit allocation indices''}, John Wiley \& Sons.
%
\bibitem{HKW}
Hemmecke, R., Koppe, M., Lee, J., Weismantel, R. (2010). ``Nonlinear integer programming",  \emph{In: Junger, M., et al. (eds.) 50 Years of Integer Programming 1958-2008: From The Early Years and State-of-the-Art}, Chapter 15, pp. 561-618, Springer, Berlin.
%
\bibitem{H}
Hochbaum, D.S. (2007).  ``Complexity and algorithms for nonlinear optimization
problems'',  \emph{Annals of Operations Research} 153, pp. 257-296.
%
\bibitem{HS}
Hochbaum, D.S., Shanthikumar, J.G. (1990).  ``Convex separable optimization is not much harder than linear optimization'',  \emph{Journal of the Association for Computing Machinery} 37, pp. 843-862.
%
\bibitem{IK}
Ibaraki, T., Katoh, N. (1988).  ``Resource Allocation Problems: Algorithmic Approaches",  \emph{MIT Press}, Boston.
%
\bibitem{IMM}
Iwata, S., Moriguchi, S., Murota, K. (2005). ``A capacity scaling algorithm for M-convex sub- modular flow",  \emph{Mathematical Programming} 103, pp.181-202.
%
\bibitem{IS}
Iwata, S., Shigeno, M. (2002). ``Conjugate scaling algorithm for Fenchel-type duality in discrete convex optimization", \emph{SIAM Journal on Optimization} 13, 204-211.
%
\bibitem{KV1987}
Katehakis, M. N., and Veinott Jr, A. F. (1987). ``The multi-armed bandit problem: decomposition and computation'', \emph{Mathematics of Operations Research}, 12(2), 262-268.

\bibitem{KSI}
Katoh, N., Shioura, A., Ibaraki, T. (2013).  ``Resource allocation problems",  \emph{In: Pardalos, P.M., Du, D.-Z., Graham, R.L. (eds.) Handbook of Combinatorial Optimization}, 2nd ed., Vol. 5, pp. 2897-2988, Springer, Berlin.
%


\bibitem{Kolesar1967}
Peter J. Kolesar (1967).``A Branch and Bound Algorithm for the Knapsack Problem'', \emph{Management Science}, 13(9), pp.723-735.

\bibitem{LL}
Lee, J., Leyffer, S. (2012), ``Mixed Integer Nonlinear Programming", \emph{The IMA Volumes in Mathematics and its Applications 154}, Springer, Berlin.
%
\bibitem{LHR}
Lewis, H. R    and Papadimitriou, C. H (1998). `` Elements of the
Theory of Computation'',  \emph{ACM Press,} New York.
%
\bibitem{MMTT}
Moriguchi, S., Murota, K., Tamura, A., and Tardella, F. (2017). ``Discrete Midpoint Convexity",  \emph{https://arxiv.org/pdf/1708.04579.pdf}.
%
\bibitem{O}
Onn, S. (2010). ``Nonlinear Discrete Optimization: An Algorithmic Theory",  \emph{European Mathematical Society}, Zurich.
%
\bibitem{PCH}
Papadimitriou, C. H. and Steiglitz, K. (1982). ``Combinatorial
Optimization: Algorithms and Complexity'', \emph{Prentice-Hall}.

%
\bibitem{PCH}
Papadimitriou , C. H. and Yannakakis, M. (1982). `` The complexity of
facets (and some facets of complexity)'', \emph{Annual ACM Symposium
on Theory of Computing archive Proceedings of the fourteenth annual
ACM symposium on Theory of computing,} San Francisco, California,
United States, pp. 255 - 260.
%
\bibitem{RossTsang1989}
Ross, K. W., and Tsang, D. H. K. (1989). ``The stochastic knapsack problem", \emph{IEEE Transactions on Communications}, 37(7), pp. 740-747.

\bibitem{ShapiroWanger1967}
Shapiro, J. F., and  H. M. Wanger, H.M. (1967). ``A finite renewal algorithm for the knapsack and turnpike models", \emph{Operations Research}, 15, pp. 319-341.
%
\bibitem{TDM}
Topkis D. M. (1978). ``Minimizing a Submodular Function on a
Lattice'',  \emph{Operations Research,}  26,  pp. 305-321.

%
\bibitem{PCH}
Yannakakis  M. (1982).  ``The complexity of the partial order dimension
problem'', \emph{SIAM J. Alg. Dis. Meth.},3,  pp. 351-370.

%
\bibitem{YM}
 Yannakakis M. (1990). ``The analysis of local search problems
modified greedy heuristic for the set covering problem and their
heuristics'',  \emph{STAG'S -90},  pp. 298-311.
%
\bibitem{Z}
Zipkin P. (2016). ``Some specially structured assemble-to-order systems",
\emph{Operations Research Letters}, 44, pp. 136-142.

%%%%%%%%%%%%the following is added referecce











\end{thebibliography}
%\section*{References}
%\bibliographystyle{amsplain}

\end{document}